\documentclass[11pt,letterpaper,reqno]{amsart}
\PassOptionsToPackage{table}{xcolor}
\usepackage{graphicx,pdfsync,amssymb,mathrsfs,amsfonts,amsbsy,color,esint,mathtools,tikz,mdframed}
\usetikzlibrary{positioning,arrows.meta,calc}
\colorlet{figVortexOne}{red!80!black}
\colorlet{figVortexTwo}{orange!90!black}
\colorlet{figVortexThree}{blue!75!black}
\colorlet{figPositive}{figVortexOne}
\colorlet{figNegative}{figVortexThree}
\colorlet{figTransition}{figVortexThree}
\colorlet{figWaiting}{blue!60!red}
\colorlet{figInitial}{black!60}
\tikzset{
  paper figure/.style={font=\small, line width=.65pt,
    line cap=round, line join=round,
    >={Stealth[length=3pt,width=3pt]}},
  fig axis/.style={draw=black!45,line width=.4pt},
  fig guide/.style={draw=black!25,line width=.4pt,dashed},
  fig geometry/.style={draw=black!55,line width=.65pt},
  fig curve/.style={line width=.9pt},
  fig selected/.style={line width=1.5pt},
  fig vector/.style={->,draw=black!55,line width=.65pt},
  fig movement/.style={->,line width=.4pt,opacity=.45,
    >={Stealth[length=2.2pt,width=2.2pt]}},
  fig vortex/.style={circle,draw=none,inner sep=0pt,
    outer sep=0pt,minimum size=4pt},
  fig particle/.style={circle,fill=black,draw=none,
    inner sep=0pt,outer sep=0pt,minimum size=2.2pt},
  fig origin/.style={circle,fill=black!55,draw=none,
    inner sep=0pt,outer sep=0pt,minimum size=1.6pt},
  fig label/.style={text=black,inner sep=2pt},
  every label/.append style={fig label}
}
\usepackage[shortlabels]{enumitem}
\usepackage{makecell}
\usepackage{booktabs,array,tabularx}

\usepackage[table]{xcolor}
\usepackage{hyperref}
\hypersetup{colorlinks=true,linkcolor=blue,filecolor=blue,urlcolor=blue,citecolor=blue}
\usepackage[utf8]{inputenc}
\usepackage[T1]{fontenc}
\usepackage{lmodern}
\usepackage{float}

\DeclareMathOperator{\dist}{dist}

\newtheorem{theorem}{Theorem}[section]
\newtheorem{lemma}[theorem]{Lemma}
\newtheorem{proposition}[theorem]{Proposition}

\setlist[enumerate]{itemsep=3pt}
\setlist[itemize]{itemsep=3pt}

\def\R{\mathbb{R}}

\numberwithin{equation}{section}
\allowdisplaybreaks[1]
\title{Large diameter growth for vortex patches}
\author{Junhong Fan}
\address{Institute of Applied Mathematics, AMSS, Chinese Academy of Sciences,
Beijing 100190; and University of Chinese Academy of Sciences,
Beijing 100049, P.R. China}
\email{fanjunhong@amss.ac.cn}

\author{Dengjun Guo}
\address{Academy of Mathematics and Systems Science,
Chinese Academy of Sciences, Beijing, China}
\email{djguo@amss.ac.cn}

\author{Xiaoyutao Luo}
\address{Academy of Mathematics and Systems Science,
Chinese Academy of Sciences, Beijing, China}
\email{xiaoyutao.luo@amss.ac.cn}
\date{}

\mathtoolsset{showonlyrefs}
\begin{document}
\begin{abstract}
For every $N>0$, the 2D incompressible Euler
equation admits a smooth, simply connected vortex patch
initially contained in the unit disk whose diameter exceeds
$N$ at some finite time. We discover a point-vortex
mechanism in which a bounded motion of three equal positive
vortices drives a passive particle outward through successive
transitions. A finite-time point-vortex approximation
transfers the resulting separation to a single simply
connected patch.
\end{abstract}
\maketitle
\section{Introduction}\label{sec:intro}

We consider the 2D incompressible Euler equation
in vorticity form on $\mathbb R^2$,
\begin{equation}\label{eq:euler}
 \partial_t\omega(t,x)+u(t,x)\cdot\nabla\omega(t,x)=0,
\end{equation}
where the velocity is recovered from the vorticity by the
Biot--Savart law
\begin{equation}\label{eq:biot-savart}
 u(t,x)=\int_{\mathbb R^2}K(x-y)\omega(t,y)\,dy.
\end{equation}
Here
\begin{equation}\label{eq:kernel}
 K(x):=\frac{x^\perp}{2\pi|x|^2},
 \qquad x^\perp:=(-x_2,x_1),\qquad x\ne0.
\end{equation}
For $r>0$ and $z\in\mathbb R^2$, we write
$B_r(z):=\{x\in\mathbb R^2:|x-z|<r\}$ and $B_r:=B_r(0)$.

For initial vorticity in $L^1(\mathbb R^2)\cap L^\infty(\mathbb R^2)$,
Yudovich's theory gives a unique global solution
\cite{Yudovich1963,MajdaBertozzi}.
In particular, initial data $\omega_0=\mathbf1_{D_0}$ remain
a vortex patch $\omega(t)=\mathbf1_{D_t}$, transported by
the area-preserving Euler flow. Chemin \cite{Chemin} and
Bertozzi--Constantin \cite{BertozziConstantin1993} proved
that an initially smooth patch boundary remains smooth
for all finite times.
Since then, vortex-patch dynamics have been studied
extensively; see, for example,
\cite{Serfati1994,Danchin1997,HmidiMateuVerdera2013,
Denisov2015,ElgindiJeong2020,HassainiaMasmoudiWheeler2020,
BertiHassainiaMasmoudi2023,KiselevLuo2023}
and the references therein.

The long-time spreading of vorticity of one sign is a
classical problem in 2D Euler dynamics.
Marchioro \cite{Marchioro1994} established an upper bound
of order $t^{1/3}$ for the growth of the vorticity support,
later sharpened to $(t\log t)^{1/4}$
\cite{IftimieSiderisGamblin1999} and, recently, $t^{1/4}$
\cite{CaoFanWang2026}.
Whether a vortex patch can have diameter tending to infinity
as $t\to\infty$ remains a longstanding open problem, even
for disconnected patches \cite{ChoiJeong2022}. We remark that for sign-changing vorticity,
\cite{IftimieSiderisGamblin1999} 
constructed solutions whose support grows linearly in time.

Here we report arbitrarily large finite-time diameter
growth: for every $N>0$, the construction yields
a smooth, simply connected patch initially contained in
the unit disk whose diameter exceeds $N$ at some finite time.

\begin{theorem}\label{thm:main}
For every $N>0$, there exist a bounded simply connected
domain $D_0\subset B_1$ with $C^\infty$ boundary and a finite
time $T_N>0$ such that the unique Euler solution
\begin{equation}\label{main:patch}
 \omega(t)=\mathbf1_{D_t},\qquad t\ge0,
\end{equation}
satisfies
\begin{equation}\label{main:diameter}
 \operatorname{diam}D_{T_N}>N.
\end{equation}
\end{theorem}

\subsection*{About the use of AI}\label{sec:ai-usage}
We used ChatGPT 6.1 sol in developing the point-vortex
construction, checking and editing the manuscript,
and preparing the figures.

The model first discovered a four-point-vortex system
producing outward motion of a passive particle.
The initial ODE proof was difficult to follow. The authors studied this construction and, through
extensive discussions with the model, discovered the
three-point-vortex mechanism and developed its present
formulation.

The authors wrote the manuscript, which was subsequently
checked and edited with ChatGPT 6.1 sol under their
guidance. The figures were also generated with AI
assistance under human instructions. The authors checked
the mathematical arguments, reviewed and approved all
text and figures, and take full responsibility for the work.



\subsection{Related work}\label{sec:related-work}

The dynamics of three and four point vortices have been
studied extensively
\cite{Aref1979,ArefPomphrey1982,Ziglin1980};
see also \cite{Newton2001,Aref2007}.
Passive-particle transport driven by three moving vortices
is particularly relevant to our construction
\cite{NeufeldTel1997,KuznetsovZaslavsky1998,
KuznetsovZaslavsky2000}.

Point-vortex dynamics have also been used to construct
Euler flows with growing support.
Zbarsky \cite{Zbarsky2021} transferred a self-similarly
expanding configuration of three point vortices of 
\emph{both signs} to Euler flows whose
vortex cores separate at rate $\sqrt{t}$.

A closely related model with circulations
$(1,1,1,\varepsilon)$ was studied near three-vortex saddle
connections in \cite{CastillaMoauroNegriniOliva1993}.
For each fixed $\varepsilon>0$, conservation of angular
impulse gives a time-independent upper bound on the
vortex positions, depending on $\varepsilon$.
Here we discover, for sufficiently small
$\varepsilon>0$, an orbit attaining a prescribed large
finite separation between the fourth vortex and the
three-vortex cluster, and transfer this separation
to a positive vortex patch.

The transfer relies on the finite-time concentration
theory of Marchioro and Pulvirenti
\cite{MarchioroPulvirenti1993,MarchioroPulvirenti1994} and its extension to
concentrated vorticity accompanied by small perturbations
in \cite{Guo}.
The latter allows us to join the concentrated vortex
cores by thin corridors, obtaining a single smooth,
simply connected patch while preserving the separation.

For further background and related work on vorticity
spreading, vortex-patch dynamics, and small-scale formation,
see \cite{Chemin2004,IftimieSiderisGamblin1999,
KiselevSverak2014,Zlatos2015,Kiselev2018,ButtaMarchioro2018,
ElgindiJeong2020,ChoiJeongPerimeter2021,ChoiJeong2022,
DrivasElgindiJeong2024,ChoiJeongWinding2025}.

\subsection{Outline of the construction}
We describe the construction in five steps.

\medskip
\textbf{Step 1: Reduction to passive-particle transport.}
We first use three vortices of unit circulation to drive
a passive particle $X(t)$ so that $|X(t)|$ grows by a
prescribed large factor, while the vortices remain bounded.
We then give the particle a small positive circulation
and obtain the required patch by the classical point-vortex
approximation of Marchioro and Pulvirenti
\cite{MarchioroPulvirenti1993} and its extension in \cite{Guo}.

\medskip
\textbf{Step 2: Three-vortex motion near the separatrix.}
For three vortices of unit circulation at positions $a_i$,
we place the center of vorticity at the origin, $\sum_i a_i=0$.
The conserved second moment and vortex energy are
\[
 I:=\frac12\sum_i|a_i|^2,\qquad
 E:=-\frac1{2\pi}\sum_{i<j}\log|a_i-a_j|.
\]
Fixing $I$ sets the size of the configuration. Removing
the overall rotation then leaves one degree of freedom,
so the motions are classified by level curves of $E$
\cite{Novikov1975,Aref1979}.

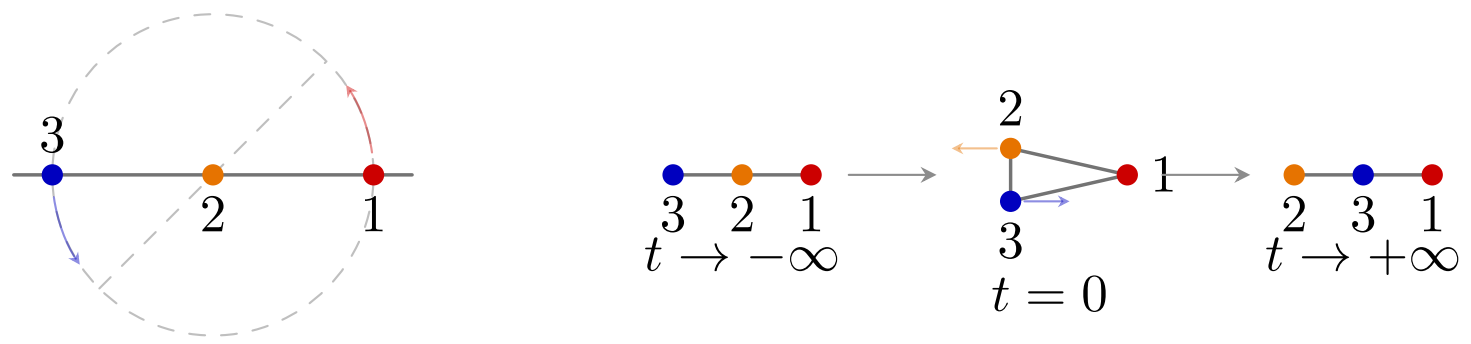
\begin{figure}[H]
\centering
\begin{tikzpicture}[paper figure]
\begin{scope}[scale=.85]
 \draw[fig guide] (0,0) circle (1.25);
 \draw[fig geometry] (-1.55,0)--(1.55,0);
 \draw[fig guide,rotate=45] (-1.25,0)--(1.25,0);
 \draw[fig movement,draw=figVortexOne] (8:1.25)
  arc[start angle=8,end angle=34,radius=1.25];
 \draw[fig movement,draw=figVortexThree] (188:1.25)
  arc[start angle=188,end angle=214,radius=1.25];
 \node[fig vortex,fill=figVortexThree,label=above:{$3$}] at (-1.25,0) {};
 \node[fig vortex,fill=figVortexTwo,label=below:{$2$}] at (0,0) {};
 \node[fig vortex,fill=figVortexOne,label=below:{$1$}] at (1.25,0) {};
\end{scope}
\begin{scope}[xshift=3.5cm,scale=.83]
 \begin{scope}
  \draw[fig geometry] (-.55,0)--(.55,0);
  \node[fig vortex,fill=figVortexThree,label=below:{$3$}] at (-.55,0) {};
  \node[fig vortex,fill=figVortexTwo,label=below:{$2$}] at (0,0) {};
  \node[fig vortex,fill=figVortexOne,label=below:{$1$}] at (.55,0) {};
  \node at (0,-.65) {$t\to-\infty$};
 \end{scope}
 \draw[fig axis,->] (.85,0)--(1.55,0);
 \begin{scope}[xshift=2.45cm]
  \draw[fig geometry] (.62,0)--(-.31,.211)--(-.31,-.211)--cycle;
  \draw[fig movement,draw=figVortexTwo] (-.42,.211)--(-.78,.211);
  \draw[fig movement,draw=figVortexThree] (-.20,-.211)--(.16,-.211);
  \node[fig vortex,fill=figVortexOne,label=right:{$1$}] at (.62,0) {};
  \node[fig vortex,fill=figVortexTwo,label=above:{$2$}] at (-.31,.211) {};
  \node[fig vortex,fill=figVortexThree,label=below:{$3$}] at (-.31,-.211) {};
  \node at (0,-.95) {$t=0$};
 \end{scope}
 \draw[fig axis,->] (3.35,0)--(4.05,0);
 \begin{scope}[xshift=4.95cm]
  \draw[fig geometry] (-.55,0)--(.55,0);
  \node[fig vortex,fill=figVortexTwo,label=below:{$2$}] at (-.55,0) {};
  \node[fig vortex,fill=figVortexThree,label=below:{$3$}] at (0,0) {};
  \node[fig vortex,fill=figVortexOne,label=below:{$1$}] at (.55,0) {};
  \node at (0,-.65) {$t\to+\infty$};
 \end{scope}
\end{scope}
\node at (0,-1.8) {(a) Collinear motion};
\node at (5.53,-1.8) {(b) Separatrix motion};
\end{tikzpicture}
\caption{(a) Collinear motion. (b) Separatrix motion in
rotating coordinates. Numbers label the vortices.}
\label{fig:intro-reference-motions}
\end{figure}

The \emph{collinear motion} is the uniform rotation of
three equally spaced vortices. The same values of $I$
and $E$ also admit a \emph{separatrix motion}: in rotating
coordinates, it approaches collinear configurations as
$t\to\pm\infty$, with vortices $2$ and $3$ exchanged;
see Figure~\ref{fig:intro-reference-motions}.

We use a family of periodic motions approaching the
separatrix. By scaling, we fix $E=-\log4/(4\pi)$ and
write $I=s^2$, calling this the \emph{$s$-vortex motion};
see Section~\ref{sec:source}. The case $s=1$ is the
separatrix. For $s>1$ close to $1$, the induced velocity
field is periodic in rotating coordinates. Its period
$T_s$ tends to infinity as $s\to1^+$, since the vortices
spend increasingly long near the collinear configuration.

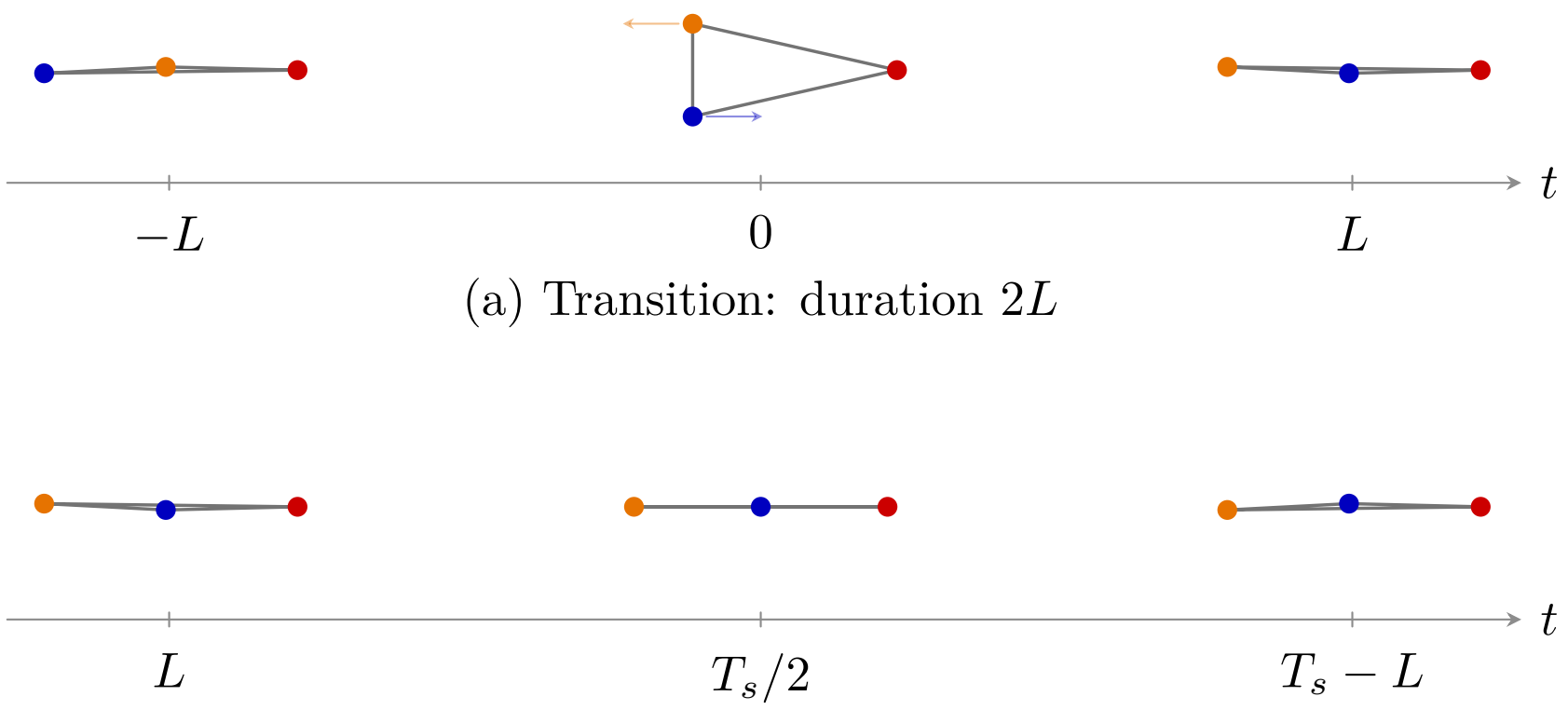
\begin{figure}[H]
\centering
\begin{tikzpicture}[paper figure]
 \node at (4.2,-1.65) {(a) Transition: duration $2L$};
 \foreach \shift/\xx/\yy/\movement in
  {0/.48/.025/0,4.2/0/.366/1,8.4/-.48/.025/0}{
  \pgfmathsetmacro{\rr}{sqrt((1-\xx*\xx-\yy*\yy)/3)}
  \begin{scope}[xshift=\shift cm,scale=.9]
   \draw[fig geometry] ({2*\rr},0)--({-\rr+\xx},\yy)
    --({-\rr-\xx},{-\yy})--cycle;
   \ifnum\movement=1
    \draw[fig movement,draw=figVortexTwo]
     ({-\rr-.11},\yy)--({-\rr-.55},\yy);
    \draw[fig movement,draw=figVortexThree]
     ({-\rr+.11},{-\yy})--({-\rr+.55},{-\yy});
   \fi
   \node[fig vortex,fill=figVortexOne] at ({2*\rr},0) {};
   \node[fig vortex,fill=figVortexTwo] at ({-\rr+\xx},\yy) {};
   \node[fig vortex,fill=figVortexThree] at ({-\rr-\xx},{-\yy}) {};
  \end{scope}
 }
 \draw[fig axis,->] (-1.15,-.8)--(9.6,-.8) node[right] {$t$};
 \foreach \pos/\txt in {0/{-L},4.2/{0},8.4/{L}}{
  \draw[fig axis] (\pos,-.75)--(\pos,-.85);
  \node[below] at (\pos,-.9) {$\txt$};
 }
 \begin{scope}[yshift=-3.1cm]
  \node at (4.2,-1.65)
   {(b) Waiting near the collinear motion: duration $T_s-2L\gg2L$};
  \foreach \shift/\xx/\yy in
   {0/-.48/.025,4.2/-.5/0,8.4/-.48/-.025}{
   \pgfmathsetmacro{\rr}{sqrt((1-\xx*\xx-\yy*\yy)/3)}
   \begin{scope}[xshift=\shift cm,scale=.9]
    \draw[fig geometry] ({2*\rr},0)--({-\rr+\xx},\yy)
     --({-\rr-\xx},{-\yy})--cycle;
    \node[fig vortex,fill=figVortexOne] at ({2*\rr},0) {};
    \node[fig vortex,fill=figVortexTwo] at ({-\rr+\xx},\yy) {};
    \node[fig vortex,fill=figVortexThree] at ({-\rr-\xx},{-\yy}) {};
   \end{scope}
  }
  \draw[fig axis,->] (-1.15,-.8)--(9.6,-.8) node[right] {$t$};
  \foreach \pos/\txt in {0/{L},4.2/{T_s/2},8.4/{T_s-L}}{
   \draw[fig axis] (\pos,-.75)--(\pos,-.85);
   \node[below] at (\pos,-.9) {$\txt$};
  }
 \end{scope}
\end{tikzpicture}
\caption{One cycle: a short transition followed by a long
wait near the collinear configuration.}
\label{fig:intro-source-cycle}
\end{figure}

Choose $1\ll L\ll T_s$. We divide each cycle into a
\emph{transition interval} $[-L,L]$, during which the
motion closely follows the separatrix, and a much longer
\emph{waiting interval} $[L,T_s-L]$, during which it
stays near the collinear motion; see Figure~\ref{fig:intro-source-cycle}.

\medskip
\textbf{Step 3: Energy gain during a transition.}
Let $Z(t)$ denote a passive particle trajectory in
rotating coordinates. To measure outward motion, we use
the time-independent particle Hamiltonian $H_*$ of the
collinear motion at the same $s$; see \eqref{ext:reference}.
For large $|Z|$, we have $H_*(Z)\approx|Z|^2$.

\begin{figure}[H]
\centering
\begin{tikzpicture}[paper figure]
 \foreach \shift/\direction/\col/\sgn/\heading in
  {0/15/figPositive/1/{Increase of $H_*$},6.2/105/figNegative/-1/{Decrease of $H_*$}}{
  \begin{scope}[xshift=\shift cm]
   \foreach \start/\shade in
    {0/figPositive,90/figNegative,180/figPositive,270/figNegative}{
    \pgfmathsetmacro{\finish}{\start+30}
    \path[fill=\shade!12]
     (\start:.85)--(\start:1.65)
     arc[start angle=\start,end angle=\finish,radius=1.65]
     --(\finish:.85)
     arc[start angle=\finish,end angle=\start,radius=.85]--cycle;
   }
   \draw[fig guide] (0,0) circle (.85);
   \draw[fig axis] (0,0)--(1.8,0);
   \draw[fig axis] (0,0)--(30:1.65);
   \draw[fig geometry] (0:.36) arc[start angle=0,end angle=30,radius=.36];
   \node at (15:.68) {$\pi/6$};
   \node[below] at (1.8,0) {$\alpha$};
   \node[fig origin,label=below left:{$0$}] at (0,0) {};
   \node[fig particle,label=above right:{$Z(-L)$}] at (\direction:1.35) {};
   \node at (0,2) {\heading};
   \begin{scope}[xshift=-1.65cm,yshift=-2.75cm]
    \fill[black!5] (0,-.65) rectangle (.7,.65);
    \draw[fig guide] (0,0)--(3.4,0);
    \draw[fig axis,->] (0,-.85)--(3.65,-.85) node[right] {$t$};
    \draw[fig axis,->] (0,-.85)--(0,.85);
    \node[above] at (0,.9) {$H_*(Z(t))$};
    \draw[fig curve,draw=\col]
     (0,0)..controls(.12,{.03*\sgn})and(.18,{.26*\sgn})
     ..(.29,{.22*\sgn})
     ..controls(.38,{.16*\sgn})and(.48,{.44*\sgn})
     ..(.7,{.45*\sgn})--(3.4,{.45*\sgn});
    \foreach \pos/\txt in {0/{-L},.7/{L},3.4/{T_s-L}}{
     \draw[fig axis] (\pos,-.80)--(\pos,-.90);
     \node[below] at (\pos,-.93) {$\txt$};
    }
    \node[below] at (.35,-1.35) {transition};
    \node[below] at (2.05,-1.35) {waiting};
   \end{scope}
  \end{scope}
 }
\end{tikzpicture}
\caption{Positive (red) and negative (blue) Hamiltonian
increments over one cycle. Vertical variations are exaggerated.}
\label{fig:intro-sectors}
\end{figure}
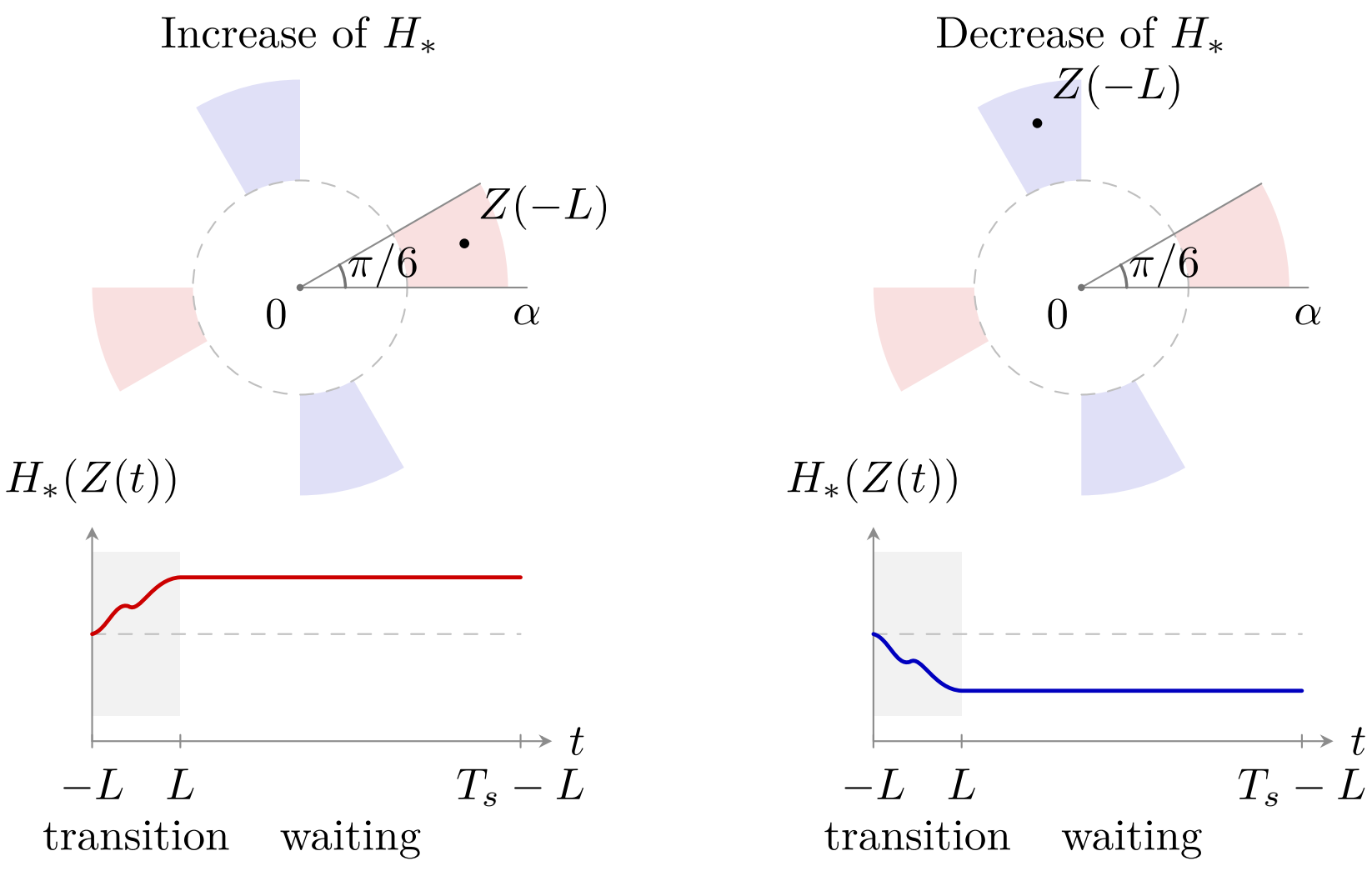

For the $s$-vortex motion, $H_*$ changes mainly during
the transition and stays almost constant during the
long wait. We find an angle-dependent gain: for $s>1$
sufficiently close to $1$ and a suitable choice of $L$,
particles starting at radius $\rho=|Z(-L)|\approx L$
in the red sectors of Figure~\ref{fig:intro-sectors} satisfy
\[
H_*(Z(T_s-L))-H_*(Z(-L))\gtrsim\rho^{-2}>0.
\]
The one-period estimate is given in Section~\ref{sec:one-period}.

\medskip
\textbf{Step 4: Winding and return during the long wait.}
We consider an arc on a level curve of $H_*$ near radius
$L$, spanning a red sector $U$. The transition gives its
endpoints unequal gains, hence different mean angular speeds.

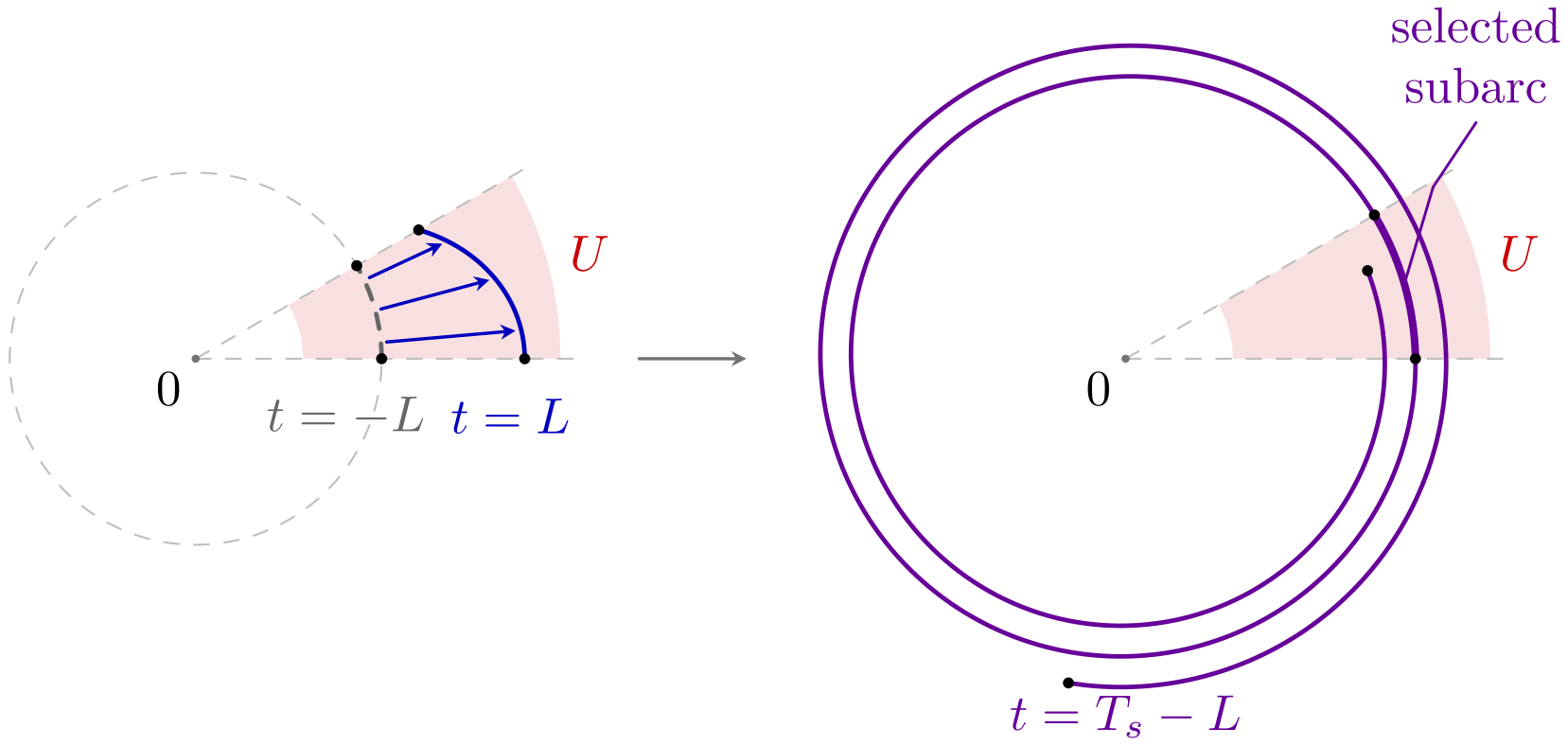
\begin{figure}[H]
\centering
\begin{tikzpicture}[paper figure]
 \begin{scope}
  \path[fill=figPositive!12]
   (0:.75)--(0:2.55)
   arc[start angle=0,end angle=30,radius=2.55]
   --(30:.75)
   arc[start angle=30,end angle=0,radius=.75]--cycle;
  \draw[fig guide] (0,0) circle (1.3);
  \draw[fig guide] (0,0)--(2.65,0);
  \draw[fig guide] (0,0)--(30:2.65);
  \node[fig origin,label=below left:{$0$}] at (0,0) {};
  \node[figPositive] at (15:2.85) {$U$};
  \draw[fig curve,draw=figInitial,dashed]
   (0:1.3) arc[start angle=0,end angle=30,radius=1.3];
  \foreach \ang in {5,15,25}{
   \pgfmathsetmacro{\rr}{1.3+cos(2*\ang)-.04}
   \draw[fig vector,draw=figTransition] (\ang:1.34)--(\ang:\rr);
  }
  \draw[fig curve,draw=figTransition,domain=0:30,samples=61,smooth]
   plot ({(1.3+cos(2*\x))*cos(\x)},
         {(1.3+cos(2*\x))*sin(\x)});
  \foreach \ang in {0,30}{
   \pgfmathsetmacro{\rr}{1.3+cos(2*\ang)}
   \node[fig particle] at (\ang:1.3) {};
   \node[fig particle] at (\ang:\rr) {};
  }
  \node[figInitial] at (1.05,-.4) {$t=-L$};
  \node[figTransition] at (2.2,-.4) {$t=L$};
  \node at (.55,-2.95) {(a) Transition};
 \end{scope}
 \draw[fig geometry,->] (3.1,0)--(3.85,0);
 \begin{scope}[xshift=6.5cm]
  \path[fill=figPositive!12]
   (0:.75)--(0:2.55)
   arc[start angle=0,end angle=30,radius=2.55]
   --(30:.75)
   arc[start angle=30,end angle=0,radius=.75]--cycle;
  \draw[fig guide] (0,0)--(2.65,0);
  \draw[fig guide] (0,0)--(30:2.65);
  \node[fig origin,label=below left:{$0$}] at (0,0) {};
  \node[figPositive] at (15:2.85) {$U$};
  \draw[fig curve,draw=figWaiting,domain=-100:740,samples=301,smooth]
   plot ({(2.3-(\x+100)/1680)*cos(\x)},
         {(2.3-(\x+100)/1680)*sin(\x)});
  \node[fig particle] at (-100:2.3) {};
  \node[fig particle] at (740:1.8) {};
  \draw[fig selected,draw=figWaiting,domain=360:390,samples=31,smooth]
   plot ({(2.3-(\x+100)/1680)*cos(\x)},
         {(2.3-(\x+100)/1680)*sin(\x)});
  \foreach \ang in {360,390}{
   \pgfmathsetmacro{\rr}{2.3-(\ang+100)/1680}
   \node[fig particle] at (\ang:\rr) {};
  }
  \draw[draw=figWaiting,line width=.55pt] (2.45,1.65)--(2.15,1.2)--(15:2.0173);
  \node[text=figWaiting,above,align=center] at (2.45,1.65)
   {selected\\subarc};
  \node[figWaiting] at (0,-2.5) {$t=T_s-L$};
  \node at (0,-2.95) {(b) Long waiting};
 \end{scope}
\end{tikzpicture}
\caption{(a) Unequal gains during the transition.
(b) Winding and selection of a returning subarc.}
\label{fig:intro-arc-evolution}
\end{figure}

During the waiting interval $[L,T_s-L]$, the arc winds
around the origin and crosses $U$ again. At $t=T_s-L$,
we select a new subarc inside $U$.
The particles on this subarc have nearly equal values
of $H_*$; see Figure~\ref{fig:intro-arc-evolution}.

\medskip
\textbf{Step 5: Iteration and construction of the vortex patch.}
Each cycle increases $H_*$ on the selected arc. After
finitely many cycles, a particle reaches the prescribed
large radius. Tracing it back through the selected arcs
gives a single particle trajectory with the required
growth. The reduction in Step~1 then yields the desired
simply connected patch; see Sections~\ref{sec:arc-selection}
and~\ref{sec:patch}.

\subsection*{Acknowledgments}

XL is supported by the National Natural Science Foundation of China   under Grant No. 12421001 and No. 12288201. DG is supported by NSFC No. 12671283. JF is supported by the National Key R\&D Program of China (2022YFA1005602) and the National Natural Science Foundation of China No. 12371212.

\section[The s-vortex motion]{The \texorpdfstring{$s$}{s}-vortex motion}\label{sec:source}
We demonstrate the existence of the \emph{$s$-vortex motion}
used in the particle construction. We establish its
periodicity in rotating coordinates and quantify its
short transitions and long waiting intervals near the
collinear motion.

\subsection{Notation and conventions}
We write $e_1=(1,0)$ and use $|\cdot|$ for the Euclidean
norm of vectors. For matrices, $M^T$ denotes the transpose
and $\|M\|$ the induced operator norm. For
$a,b,\xi\in\mathbb R^2$, we set
$(a\otimes b)\xi=a(b\cdot\xi)$.

The constants $c,C>0$ may change from line to line and
are universal unless their dependence is specified.
For nonnegative $A,B$, we write $A\lesssim B$ if
$A\le CB$, $A\gtrsim B$ if $B\lesssim A$, and
$A\approx B$ if both inequalities hold.
A subscript, as in $\lesssim_M$, indicates allowed
dependence on that parameter. For $G\ge0$, we write
$F=O(G)$ when $|F|\lesssim G$.

\subsection{Construction}
Consider three point vortices
$a_i(t)=(a_{i,1}(t),a_{i,2}(t))$, each with circulation one:
\begin{equation}\label{src:cartesian}
 \dot a_i(t)=\sum_{\substack{1\le j\le3\\j\ne i}}K(a_i(t)-a_j(t)),
 \qquad i=1,2,3.
\end{equation}

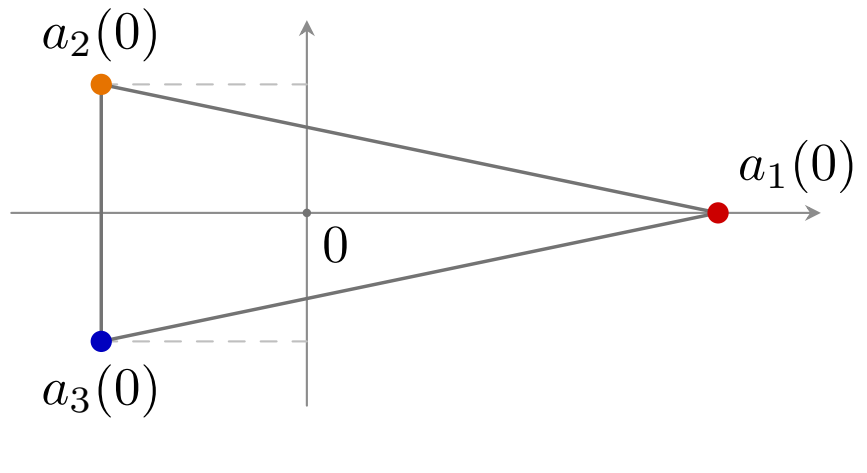
\begin{figure}[H]
\centering
\begin{tikzpicture}[paper figure,scale=.85]
 \coordinate (A1) at (3.2,0);
 \coordinate (A2) at (-1.6,1);
 \coordinate (A3) at (-1.6,-1);
 \draw[fig axis,->] (-2.3,0)--(4,0);
 \draw[fig axis,->] (0,-1.5)--(0,1.5);
 \draw[fig guide] (A2)--(0,1);
 \draw[fig guide] (A3)--(0,-1);
 \draw[fig geometry] (A1)--(A2)--(A3)--cycle;
 \node[fig origin,label=below right:{$0$}] at (0,0) {};
 \node[fig vortex,fill=figVortexOne,label=above right:{$a_1(0)$}] at (A1) {};
 \node[fig vortex,fill=figVortexTwo,label=above:{$a_2(0)$}] at (A2) {};
 \node[fig vortex,fill=figVortexThree,label=below:{$a_3(0)$}] at (A3) {};
\end{tikzpicture}
\caption{Initial positions, with
$a_{1,1}(0)=-2a_{2,1}(0)>0$.}
\label{fig:initial-positions}
\end{figure}

For each $s\ge1$, the following lemma specifies the
\emph{$s$-vortex motion}, with conserved second moment
$s^2$ and initial configuration illustrated in
Figure~\ref{fig:initial-positions}.
Throughout the paper, we usually suppress the dependence on $s$.

\begin{lemma}[Construction of $s$-vortex motion]\label{lem:vortex-invariants}
\label{lem:initial-data}
For any $s\ge1$, there exists a global solution
$
 a_i(t)=\bigl(a_{i,1}(t),a_{i,2}(t)\bigr)
$
of \eqref{src:cartesian} satisfying the following properties.
\begin{enumerate}
\item[(i)] For any $t\in\mathbb R$, there holds
\begin{equation}\label{src:invariants}
 \qquad \quad\,\,\,\, \sum_{i=1}^3a_i(t)=0,\quad
 \frac12\sum_{i=1}^3|a_i(t)|^2=s^2,\quad
 \prod_{1\le i<j\le3}|a_i(t)-a_j(t)|^2=4
\end{equation}
and
\begin{equation}\label{src:distance-bounds}
 \qquad \quad \,\,\,\,\, a_1(t)\ne0,\quad
 |a_i(t)|\le\sqrt2s,\quad
 |a_i(t)-a_j(t)|^2\ge\frac4{(6s^2)^2}
 \quad(i\ne j).
\end{equation}

\item[(ii)] The initial data satisfy
\[
 \qquad a_{1,1}(0)>0,\quad
 a_{2,1}(0)=a_{3,1}(0),\quad
 0<a_{2,2}(0)=-a_{3,2}(0)<s/2.
\]
\end{enumerate}
\end{lemma}
\begin{proof}
For distinct initial positions, \eqref{src:cartesian} has a
unique global solution because all circulations are positive,
see \cite[Proposition 3.2]{DonatiGodardCadillacIftimie}.
We observe that the center of vorticity, second moment,
and Hamiltonian are conserved:
\begin{equation}\label{eq differential conservation}
 \frac d{dt}\sum_i a_i(t)=0,\quad
 \frac d{dt}\frac12\sum_i|a_i(t)|^2=0,\quad
 \frac d{dt}\sum_{i<j}\log|a_i(t)-a_j(t)|^2=0.
\end{equation}

By Lemma~\ref{lem:shape-function-properties}\textup{(ii)},
there exists a unique $h\in(0,s/2)$ such that
\[
 h^2(3s^2-2h^2)^2=1.
\]
Set $r_0=\sqrt{(s^2-h^2)/3}$ and choose
\[
 a_1(0)=(2r_0,0),\qquad
 a_2(0)=(-r_0,h),\qquad
 a_3(0)=(-r_0,-h).
\]
These initial positions satisfy \textup{(ii)} and
\eqref{src:invariants} at $t=0$.
Conservation gives \eqref{src:invariants} for all
$t\in\mathbb R$. In particular,
$|a_i(t)|\le\sqrt2s$, and
\[
 \sum_{i<j}|a_i(t)-a_j(t)|^2
 =3\sum_i|a_i(t)|^2
  -\left|\sum_i a_i(t)\right|^2
 =6s^2.
\]
For $\{i,j,k\}=\{1,2,3\}$, 
\eqref{src:invariants} also gives
\[
 |a_i(t)-a_j(t)|^2
 =\frac4{|a_i(t)-a_k(t)|^2
          |a_j(t)-a_k(t)|^2}
 \ge\frac4{(6s^2)^2}.
\]
It remains to check that $a_1(t)\ne0$.
If $a_1(t_0)=0$ for some $t_0\in\mathbb R$, then
\eqref{src:invariants} gives
\[
 a_2(t_0)=-a_3(t_0)\ne0.
\]
These positions provide a uniformly rotating collinear solution
with $a_1(t)\equiv0$. Uniqueness for \eqref{src:cartesian}
then contradicts \textup{(ii)}.
\end{proof}
\subsection{Rotating coordinates}
\begin{figure}[H]
\centering
\begin{tikzpicture}[paper figure,x=1.65cm,y=1.35cm]
 \coordinate (O) at (0,0);
 \coordinate (M) at (-1,0);
 \coordinate (B1) at (2,0);
 \coordinate (B2) at (-.42,1);
 \coordinate (B3) at (-1.58,-1);
 \draw[fig axis,->] (-2.05,0)--(2.55,0);
 \draw[fig geometry] (B3)--(B2);
 \draw[fig vector,shorten >=3pt] (M)--(B2)
      node[midway,above left,inner sep=3pt] {$(x,y)$};
 \draw[fig axis] (-1,-.06)--(-1,.06);
 \node[fig origin,label=below right:{$0$}] at (O) {};
 \node[fig vortex,fill=figVortexOne,label=above:{$b_1$}] at (B1) {};
 \node[fig vortex,fill=figVortexTwo,label=above right:{$b_2$}] at (B2) {};
 \node[fig vortex,fill=figVortexThree,label=below left:{$b_3$}] at (B3) {};
 \draw[fig geometry,<->] (-1,-.3)--(0,-.3)
      node[midway,below=2pt] {$r$};
\end{tikzpicture}
\caption{The shape variables in rotating coordinates.}
\label{fig:source-geometry}
\end{figure}
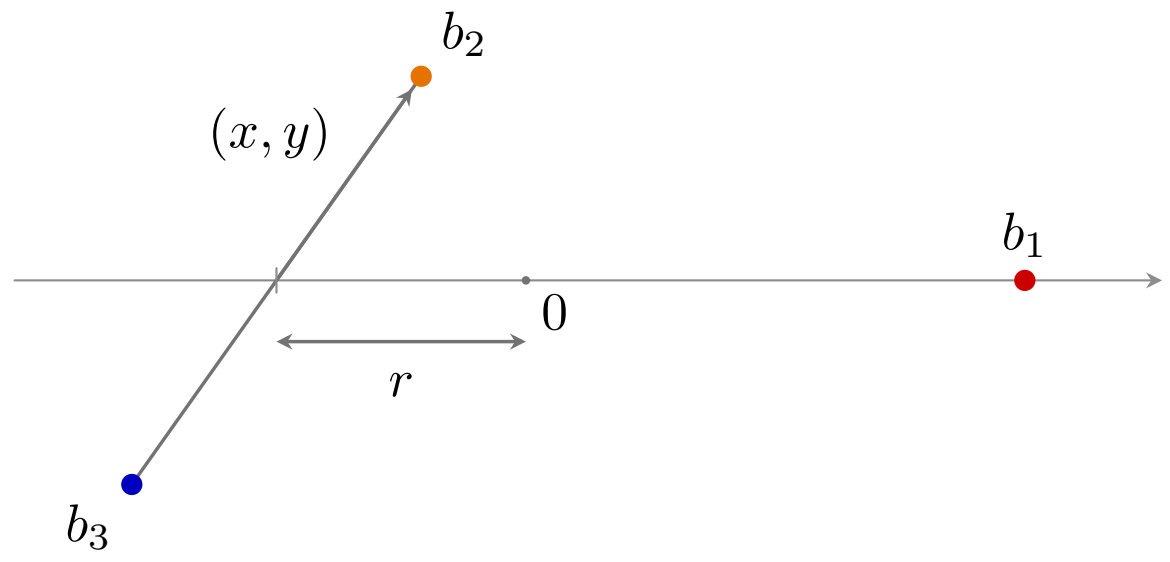

In rotating coordinates, the first vortex stays on the
positive horizontal axis; see Figure~\ref{fig:source-geometry}.
Since $a_1(t)\ne0$ by Lemma~\ref{lem:vortex-invariants},
we define
\begin{equation}\label{src:frame}
 \mathcal R(t)=\frac1{|a_1(t)|}
 \begin{pmatrix}a_{1,1}(t)&a_{1,2}(t)\\-a_{1,2}(t)&a_{1,1}(t)\end{pmatrix}    \quad \text{and} \quad b_i(t)=\mathcal R(t)a_i(t).
\end{equation}
Define
$$r(t)=\frac{|a_1(t)|}{2}>0 \qquad \text{and} \qquad (x(t),y(t))=\frac{b_2(t)-b_3(t)}2,$$
then
\eqref{src:frame} and \eqref{src:invariants} give
\begin{equation}\label{src:geometry}
 b_1(t)=(2r(t),0), \,\,\,\, b_2(t)=(-r(t)+x(t),y(t)), \,\,\,\, b_3(t)=(-r(t)-x(t),-y(t)).
\end{equation}

The reduced motion in $(x,y)$ is Hamiltonian, with
Hamiltonian $F_s$; see \eqref{src:shape-ode} below.
For $s\ge1$ and $(x,y)\in\mathbb R^2$, define
\begin{equation}\label{src:shape-function}
 F_s(x,y):=(x^2+y^2)\bigl[(3s^2-2x^2-2y^2)^2
                  -12(s^2-x^2-y^2)x^2\bigr].
\end{equation}
The properties of $F_s$ used below are collected in
Lemma~\ref{lem:shape-function-properties}.

\begin{lemma}\label{lem:shape-equations}
For any $s\ge1$, the s-vortex motion $a_i(t)$ in
Lemma~\ref{lem:vortex-invariants} satisfies the following properties
for every $t\in\mathbb R$.
\begin{enumerate}
\item[(i)] The coordinates defined above  satisfy $F_s(x(t),y(t))=1$ and
\begin{equation}\label{src:size}
 3r(t)^2+x(t)^2+y(t)^2=s^2.
\end{equation}
\item[(ii)] The rotation matrix $\mathcal R(t)$ satisfies
\begin{equation}\label{src:frame-derivative}
 \dot{\mathcal R}(t)\mathcal R^T(t)
 =-\frac{a_1(t)^\perp\cdot\dot a_1(t)}{|a_1(t)|^2}
 \begin{pmatrix}0&-1\\1&0\end{pmatrix},
\end{equation}
where
\begin{equation}\label{src:rotation-speed}
 \begin{aligned}
 \frac{a_1(t)^\perp\cdot\dot a_1(t)}{|a_1(t)|^2}
 =&\frac1{4\pi r(t)}
 \Biggl[
 \frac{3r(t)-x(t)}{(3r(t)-x(t))^2+y(t)^2}\\
 &\hphantom{\frac1{4\pi r(t)}\Biggl[}
 {}+\frac{3r(t)+x(t)}{(3r(t)+x(t))^2+y(t)^2}
 \Biggr].
 \end{aligned}
\end{equation}

\item[(iii)] The equations for $x(t)$ and $y(t)$ are
\begin{equation}\label{src:shape-ode}
 \left(\dot x(t),\dot y(t)\right)
 =\frac1{8\pi}
 \bigl(-\partial_yF_s(x(t),y(t)),
        \partial_xF_s(x(t),y(t))\bigr).
\end{equation}
\end{enumerate}
\end{lemma}
\begin{proof}
Equations \eqref{src:invariants}, \eqref{src:frame} and
\eqref{src:geometry} give
$$
 2s^2=\sum_i|b_i(t)|^2
 =6r(t)^2+2x(t)^2+2y(t)^2,
$$
which proves \eqref{src:size}.
The orthogonality of $\mathcal R(t)$,
\eqref{src:invariants} and \eqref{src:geometry} also give
\begin{equation}\label{eq 1 normalize}
 \begin{aligned}
 1&=\frac14\prod_{i<j}|b_i(t)-b_j(t)|^2\\
  &=(x(t)^2+y(t)^2)
    \bigl[(3r(t)-x(t))^2+y(t)^2\bigr]
    \bigl[(3r(t)+x(t))^2+y(t)^2\bigr]\\
  &=(x(t)^2+y(t)^2)
    \bigl[(9r(t)^2+x(t)^2+y(t)^2)^2
          -36r(t)^2x(t)^2\bigr].
 \end{aligned}
\end{equation}
Substituting \eqref{src:size} into \eqref{eq 1 normalize}
and using \eqref{src:shape-function} gives
$F_s(x(t),y(t))=1$, as stated in \textup{(i)}.

Differentiating \eqref{src:frame} gives
\eqref{src:frame-derivative}.
Since $\mathcal R(t)K(x)=K(\mathcal R(t)x)$ for $x\ne0$,
equations \eqref{src:cartesian}, \eqref{src:frame} and
\eqref{src:frame-derivative} imply
\begin{equation}\label{src:rotating}
 \begin{split}
 \dot b_i(t)
 &=\dot{\mathcal R}(t)a_i(t)
   +\mathcal R(t)\dot a_i(t)\\
 &=\sum_{j\ne i}K(b_i(t)-b_j(t))
   -\frac{a_1(t)^\perp\cdot\dot a_1(t)}
          {|a_1(t)|^2}b_i(t)^\perp.
 \end{split}
\end{equation}
Taking the second component of \eqref{src:rotating} for $i=1$
and using \eqref{src:geometry}, we obtain
\[
 \begin{aligned}
 2r(t)\frac{a_1(t)^\perp\cdot\dot a_1(t)}{|a_1(t)|^2}={}&\frac1{2\pi}\left(
 \frac{3r(t)-x(t)}{(3r(t)-x(t))^2+y(t)^2}
 +\frac{3r(t)+x(t)}{(3r(t)+x(t))^2+y(t)^2}\right).
 \end{aligned}
\]
Since $r(t)>0$, this gives \eqref{src:rotation-speed}
and hence \textup{(ii)}.
By \eqref{src:geometry},
\[
 (\dot x(t),\dot y(t))
 =\frac12\bigl(\dot b_2(t)-\dot b_3(t)\bigr).
\]
Using \eqref{src:size}, \eqref{src:rotation-speed}, and
\eqref{src:rotating}, and then simplifying with
\eqref{eq 1 normalize}, we obtain
\[
 \begin{aligned}
 \dot x(t)
 =&\frac{y(t)}{4\pi}\Biggl[
 \frac{2-x(t)/r(t)}{(3r(t)-x(t))^2+y(t)^2}
 +\frac{2+x(t)/r(t)}{(3r(t)+x(t))^2+y(t)^2}
 -\frac1{x(t)^2+y(t)^2}\Biggr]\\
 ={}&-\frac{3y(t)}{4\pi}
 (s^2-2x(t)^2-2y(t)^2)
 (3s^2-6x(t)^2-2y(t)^2)
 \end{aligned}
\]
and
\[
 \begin{aligned}
 \dot y(t)
 =&\frac1{4\pi}\Biggl[
 \frac{x(t)^2/r(t)-2x(t)-3r(t)}
      {(3r(t)-x(t))^2+y(t)^2}
 +\frac{3r(t)-2x(t)-x(t)^2/r(t)}
       {(3r(t)+x(t))^2+y(t)^2}
 +\frac{x(t)}{x(t)^2+y(t)^2}\Biggr]\\
 ={}&\frac{3x(t)}{4\pi}\Bigl[
 (3s^2-4x(t)^2-4y(t)^2)
 (s^2-4x(t)^2-4y(t)^2)\\
 &\qquad \qquad
 +4(s^2-2x(t)^2-2y(t)^2)y(t)^2\Bigr].
 \end{aligned}
\]
Comparing these expressions with \eqref{src:shape-derivative-y}
and \eqref{src:shape-derivatives} yields \eqref{src:shape-ode}.
\end{proof}

\subsection{Periodicity and exponential estimates}
The cycle length $T_s$ in \eqref{src:periodicity} below
corresponds to exchanging vortices $2$ and $3$ in the
rotating frame, and each vortex returns to its position
after $2T_s$.

\begin{proposition}\label{prop:source-family}
There exists a universal constant $c>0$
such that for any $1<s<2$, the following hold for the s-vortex motion
$a_i(t)$ in Lemma~\ref{lem:vortex-invariants}.
\begin{enumerate}
\item[(i)] For any $t\in\mathbb R$,
we have $r(t)\ge s/2$.

\item[(ii)] There exists $T_s>0$ such that
$\lim_{s\to1^+}T_s=\infty$ and for every $t\in\mathbb R$,
\begin{equation}\label{src:periodicity}
\qquad \quad \,\,\, b_1(t+T_s)=b_1(t),\quad
 b_2(t+T_s)=b_3(t),\quad
 b_3(t+T_s)=b_2(t)
\end{equation}
and
\begin{equation}\label{src:periodicity-two}
 \frac{a_1(t+T_s)^\perp\cdot\dot a_1(t+T_s)}
      {|a_1(t+T_s)|^2}
 =\frac{a_1(t)^\perp\cdot\dot a_1(t)}
       {|a_1(t)|^2}.
\end{equation}

\item[(iii)] We have
\begin{equation}\label{src:source-decay}
\begin{aligned}
 |b_1(t)-se_1|+|b_2(t)+se_1|+|b_3(t)|
 &\le c^{-1}e^{-c|t|},
 &&0\le t\le T_s/2,\\
 |b_1(t)-se_1|+|b_2(t)|+|b_3(t)+se_1|
 &\le c^{-1}e^{-c|t|},
 &&-T_s/2\le t<0.
\end{aligned}
\end{equation}
Moreover, for $|t|\le T_s/2$,
\begin{equation}\label{src:speed-decay}
 \left|
 \frac{a_1(t)^\perp\cdot\dot a_1(t)}
      {|a_1(t)|^2}
 -\frac3{4\pi s^2}
 \right|
 \le\frac1c e^{-c|t|}.
\end{equation}
\end{enumerate}
\end{proposition}
\begin{proof}
\noindent\textbf{Step 1: The level curve of $F_s(x,y)$.}
Lemma~\ref{lem:vortex-invariants} and \eqref{src:frame} give
\begin{equation}\label{src:initial-x}
 x(0)=0,\qquad 0<y(0)=a_{2,2}(0)<s/2.
\end{equation}
Thus, Lemma~\ref{lem:shape-equations}\textup{(i)} and
Lemma~\ref{lem:shape-function-properties}\textup{(ii)} imply that
$y(0)=y_s(0)$.
For $s>1$, Lemma~\ref{lem:shape-function-properties}\textup{(iii)}
gives the smooth closed curve \eqref{def level curve}
containing $(x(0),y(0))$.
The vector field in \eqref{src:shape-ode} does not vanish on this curve.
Hence the solution traverses it in finite time, so
$x(t)$ and $y(t)$ are periodic and
\begin{equation}\label{eq xs ys bound s/2}
 x(t)^2+y(t)^2<s^2/4,\qquad t\in\mathbb R,
\end{equation}
which, together with \eqref{src:size} yields \textup{(i)}.

\medskip
\noindent\textbf{Step 2: Periodicity.}
For $1<s<2$, let $2T_s$ be the least positive period of
$(x(t),y(t))$.
By \eqref{src:initial-x} and Lemma~\ref{lem:shape-function-properties}\textup{(i)},
 $(-x(-t),y(-t))$ solves \eqref{src:shape-ode}
with the same initial data as $(x(t),y(t))$.
Thus,
\begin{equation}\label{src:symmetry}
 x(-t)=-x(t),\qquad y(-t)=y(t),
 \qquad t\in\mathbb R.
\end{equation}
By $2T_s$-periodicity,
\[
 x(T_s)=x(-T_s)=-x(T_s)
 \quad\Longrightarrow\quad x(T_s)=0=x(0).
\]
Using Lemma~\ref{lem:shape-equations}\textup{(i)}, we have $F_s(0,y(T_s))=1$.
Thus, Lemma~\ref{lem:shape-function-properties}\textup{(ii)} and the minimality of $2T_s$ yield $y(T_s)=-y(0)$, which gives
\begin{equation}\label{eq half period Ts}
 (x(T_s),y(T_s))
 =(0,-y(0))=(-x(0),-y(0)).
\end{equation}
Since $F_s(-x,-y)=F_s(x,y)$, both
$(x(t+T_s),y(t+T_s))$ and $(-x(t),-y(t))$
solve \eqref{src:shape-ode} with the same initial data.
Thus,
\begin{equation}\label{src:half-period-symmetry}
 x(t+T_s)=-x(t),\qquad y(t+T_s)=-y(t),
 \qquad t\in\mathbb R.
\end{equation}
Since $\dot x(0)=-\frac1{8\pi}\partial_yF_s(0,y(0))<0$
by \eqref{src:shape-ode}, \eqref{src:shape-derivative-y}
and \eqref{src:initial-x}, the uniqueness in
Lemma~\ref{lem:shape-function-properties}\textup{(ii)}
and minimality of $2T_s$ yield
\begin{equation}\label{eq xs on 0 Ts}
 x(t)<0 \quad \text{for} \quad 0<t<T_s
 \qquad \text{and} \qquad
 x(t)>0 \quad \text{for} \quad T_s<t<2T_s.
\end{equation}
Taking $t=-T_s/2$ in \eqref{src:half-period-symmetry}
and using \eqref{src:symmetry}, we obtain
\begin{equation}\label{eq ys T/2}
 y(T_s/2)=-y(-T_s/2)=-y(T_s/2)
 \quad\Longrightarrow\quad y(T_s/2)=0.
\end{equation}

Using \eqref{src:size} and \eqref{src:half-period-symmetry},
we obtain $r(t+T_s)=r(t)$, which, together with
\eqref{src:geometry} and \eqref{src:rotation-speed}, yields
\eqref{src:periodicity} and \eqref{src:periodicity-two}.

\medskip
\noindent\textbf{Step 3: The period and exponential estimates.}
All implicit constants below are uniform for $1<s<2$.
By \eqref{eq xs on 0 Ts}, \eqref{src:shape-ode} and
Lemma~\ref{lem:shape-function-properties}\textup{(iv)},
$\dot y(t)<0$ for $0<t<T_s$.
Together with \eqref{eq ys T/2}, this gives
\[
 x(t)\le0\le y(t),\qquad 0\le t\le T_s/2.
\]
For $0\le t\le T_s/2$, \eqref{src:shape-ode} and
Lemma~\ref{lem:shape-function-properties}\textup{(iv)}
give
\begin{equation}\label{eq ode for ys}
 \begin{aligned}
 -\frac d{dt}\bigl(s/2+x(t)+y(t)\bigr)
 &=\frac{\partial_yF_s(x(t),y(t))
          -\partial_xF_s(x(t),y(t))}{8\pi}\\
 &\approx y(t)+|x(t)|\bigl(y(t)+\sqrt{s-1}\bigr).
 \end{aligned}
\end{equation}
Since $|x(t)|<s/2<1$ and $y(t)\ge0$, we have
\[
 \begin{aligned}
 y(t)+|x(t)|\sqrt{s-1}
 \le y(t)+|x(t)|\bigl(y(t)+\sqrt{s-1}\bigr)\lesssim y(t)+\sqrt{s-1}.
 \end{aligned}
\]
Moreover, \eqref{src:radial-bounds} and $\sqrt{s-1}<1$ yield
\[
 \begin{aligned}
 y(t)+|x(t)|\sqrt{s-1}
 \ge\frac{\sqrt{s-1}}{2}\bigl(y(t)+|x(t)|\bigr)+\frac{y(t)}{2}\gtrsim y(t)+\sqrt{s-1}.
 \end{aligned}
\]
Therefore, since $x(t)\le0\le y(t)$, it follows from \eqref{eq ode for ys} and \eqref{fs:level-displacement} that
\[
 -\frac d{dt}\bigl(s/2+x(t)+y(t)\bigr)
 \approx s/2+x(t)+y(t),
 \qquad 0\le t\le T_s/2.
\]
Gr\"onwall's inequality then yields a constant $c>0$ independent of $s$, such that
\begin{equation}\label{src:decay-comparison}
 \bigl(s/2+y(0)\bigr)e^{-t/c}
 \le s/2+x(t)+y(t)
 \le\bigl(s/2+y(0)\bigr)e^{-ct},
\end{equation}
for $0\le t\le T_s/2$.

Since $y(T_s/2)=0$, applying \eqref{fs:level-displacement} to $(x,y)=\left( x(T_s/2),y(T_s/2) \right)$, we obtain
\[
 s/2+x(T_s/2)\approx\sqrt{s-1}.
\]
Thus, \eqref{src:decay-comparison} yields
$$
\sqrt{s-1} \gtrsim e^{-\frac{T_s}{2c}}.
$$
Hence $T_s\to\infty$ as $s\to1^+$.

Using \eqref{src:symmetry} and the upper bound in \eqref{src:decay-comparison}, we also get
\begin{equation}\label{src:shape-decay}
 s/2-|x(t)|+|y(t)|
 \lesssim e^{-c|t|},
 \qquad |t|\le T_s/2.
\end{equation}

\medskip
\noindent\textbf{Step 4: Positions and angular speed.}
For $|t|\le T_s/2$, \eqref{src:size} and \textup{(i)} give
\[
 \begin{aligned}
 0\le r(t)-s/2
 =\frac{s^2/4-x(t)^2-y(t)^2}{3(r(t)+s/2)}&\lesssim(s/2+|x(t)|)(s/2-|x(t)|)\\
 &\lesssim s/2-|x(t)|.
 \end{aligned}
\]
Combining this with \eqref{src:shape-decay}, we obtain
\[
 |r(t)-s/2|+\bigl||x(t)|-s/2\bigr|+|y(t)|
 \lesssim e^{-c|t|}.
\]
By \eqref{src:geometry} and the triangle inequality, if $x(t)\ge0$, then
$$
\begin{aligned}
&|b_1(t)-se_1|+|b_2(t)|+|b_3(t)+se_1| \\
&\lesssim |r(t)-s/2|+|x(t)-r(t)|+|s-r(t)-x(t)|+|y(t)| \\
&\lesssim |r(t)-s/2|+\bigl||x(t)|-s/2\bigr|+|y(t)|
 \lesssim e^{-c|t|}.
\end{aligned}
$$
If $x(t)<0$, the same calculation gives
\[
 \begin{aligned}
 |b_1(t)-se_1|+|b_2(t)+se_1|+|b_3(t)|\lesssim
  e^{-c|t|}.
 \end{aligned}
\]
After decreasing $c>0$ if necessary,
\eqref{src:source-decay} follows.

For the angular-speed estimate \eqref{src:speed-decay}, set
\[
 J(x,y,r)=\frac1{4\pi r}
 \left(
 \frac{3r-x}{(3r-x)^2+y^2}
 +\frac{3r+x}{(3r+x)^2+y^2}
 \right).
\]
Since $J(-x,y,r)=J(x,y,r)$ and
$J(s/2,0,s/2)=3/(4\pi s^2)$,
\eqref{src:rotation-speed} gives
\[
 \begin{aligned}
 \left|
 \frac{a_1(t)^\perp\cdot\dot a_1(t)}
      {|a_1(t)|^2}
 -\frac3{4\pi s^2}
 \right|
 &=\left|J(|x(t)|,y(t),r(t))-J(s/2,0,s/2)\right|\\
 &\le
 \left(|r(t)-s/2|+\bigl||x(t)|-s/2\bigr|+|y(t)|\right)
 \left\|\nabla J\right\|_{L^\infty(\mathcal D)},
 \end{aligned}
\]
where (using \eqref{eq xs ys bound s/2}, \eqref{src:size} and \textup{(i)})
\[
 \mathcal D
 :=\left\{(x,y,r)\in\mathbb R^3:
 |x|\le s/2,\ |y|\le s/2,\ s/2\le r\le s/\sqrt3\right\}.
\]
The denominators in $J$ are uniformly bounded away
from zero on $\mathcal D$, so
$\|\nabla J\|_{L^\infty(\mathcal D)}\lesssim1$.
The estimate follows from \eqref{src:shape-decay}.
\end{proof}
\section{Particle motion and the collinear Hamiltonian}\label{sec:coordinates}
We now study a passive particle driven by the $s$-vortex
motion. In rotating coordinates, the collinear Hamiltonian
provides a time-independent reference for measuring
outward motion.
Throughout Sections~\ref{sec:coordinates}--\ref{sec:arc-selection},
$c,C>0$ are independent of $s\in(1,2)$.

\subsection{The particle equations and the two Hamiltonians}
Fix $1<s<2$ and let $a_i(t)$ be the $s$-vortex motion
from Lemma~\ref{lem:vortex-invariants}.
A passive particle trajectory $X(t)\in\mathbb R^2$ satisfies
\begin{equation}\label{ext:particle}
 \dot X(t)=\sum_{i=1}^3K(X(t)-a_i(t)).
\end{equation}
On any interval where $X(t)\ne0$, write its coordinates
in the rotating frame \eqref{src:frame} as
\begin{equation}\label{ext:polar}
 Z(t):=\mathcal R(t)X(t)
 =\rho(t)\bigl(\cos\theta(t),\sin\theta(t)\bigr).
\end{equation}
Here $\rho(t)=|X(t)|$, and $\theta(t)$ is measured
counterclockwise, initially in $[0,2\pi)$, with full
turns counted.
The time-dependent particle Hamiltonian in this frame is
\begin{equation}\label{ext:actual-hamiltonian}
 H(t,Z)=-\frac1{2\pi}\sum_{i=1}^3\log|Z-b_i(t)|
 +\frac{|Z|^2}{2}
 \frac{a_1(t)^\perp\cdot\dot a_1(t)}{|a_1(t)|^2}.
\end{equation}
Indeed, differentiating \eqref{ext:polar} and using
\eqref{ext:particle} and \eqref{src:frame-derivative} gives
\[
 \dot Z(t)=\sum_iK(Z(t)-b_i(t))
 -\frac{a_1(t)^\perp\cdot\dot a_1(t)}{|a_1(t)|^2}
 Z(t)^\perp=-\nabla_Z^\perp H(t,Z(t)).
\]
Writing $H(t,\rho,\theta)=H(t,\rho\cos\theta,\rho\sin\theta)$,
we obtain
\begin{equation}\label{ext:polar-ode}
 \dot\rho(t)=\frac{\partial_\theta H(t,\rho(t),\theta(t))}{\rho(t)}
 \qquad \text{and} \qquad
 \dot\theta(t)=-\frac{\partial_\rho H(t,\rho(t),\theta(t))}{\rho(t)}.
\end{equation}

For the collinear motion at the same $s$, the rotating
frame has angular velocity $3/(4\pi s^2)$ and fixed
vortex positions $se_1$, $0$, and $-se_1$.
Its time-independent particle Hamiltonian is
\begin{equation}\label{ext:reference}
 H_*(Z)=-\frac{\log|Z-se_1|+\log|Z|+\log|Z+se_1|}{2\pi}
 +\frac{3|Z|^2}{8\pi s^2}
\end{equation}
We use the same notation for $H_*$ in polar coordinates.
Along the actual particle trajectory $Z(t)$,
\eqref{ext:polar-ode} gives
\begin{equation}\label{ext:energy-identity}
 \frac d{dt}H_*(Z(t))
 =\left.
 \frac{\partial_\rho H_*\,\partial_\theta H
       -\partial_\theta H_*\,\partial_\rho H}{\rho}
 \right|_{(\rho,\theta)=(\rho(t),\theta(t))}.
\end{equation}

\subsection{Hamiltonian estimates}
The following two lemmas control the radius and angular
speed on the exterior level curves of $H_*$.
The standard proofs are given in
Appendix~\ref{app:hamiltonian-estimates}.
\begin{lemma}[The exterior Hamiltonian levels]\label{lem:reference-curves}
There exist universal constants $\rho_0>0$ and $h_0>0$
such that for any $1<s<2$, there hold
\begin{enumerate}
\item[(i)] For any $\rho\ge\rho_0$ and $\theta\in\mathbb R$, we have $H_*(\rho,\theta)\approx\rho^2$,
\begin{equation}\label{ext:reference-bounds}
 \qquad\partial_\rho H_*
 =\frac{3\rho}{4\pi s^2}-\frac3{2\pi\rho}+O(\rho^{-3}),
 \quad
 \partial_\theta H_*=-\frac{s^2}{\pi\rho^2}\sin2\theta+O(\rho^{-4})
\end{equation}
and
\begin{equation}\label{ext:reference-twist}
 \partial_\rho^2H_*-\frac{\partial_\rho H_*}{\rho}
 =\frac3{\pi\rho^2}+O(\rho^{-4}).
\end{equation}

\item[(ii)] For any $\rho,\widetilde\rho\ge\rho_0$
and $\theta,\widetilde\theta\in\mathbb R$,
\begin{equation}\label{ext:energy-radius}
 |\rho-\widetilde\rho|
 \lesssim \frac{|H_*(\rho,\theta)-H_*(\widetilde\rho,\widetilde\theta)|}
              {\min(\rho,\widetilde\rho)}
 +\frac{1}{\min(\rho,\widetilde\rho)^3}.
\end{equation}

\item[(iii)] For any $h\ge h_0$ and $\theta\in\mathbb R$,
there exists a unique $\bar\rho(h,\theta)\ge\rho_0$ satisfying
\begin{equation}\label{mean:level-radius}
 H_*(\bar\rho(h,\theta),\theta)=h.
\end{equation}
This solution $\bar \rho$ depends smoothly on $(s,h,\theta)$.
\end{enumerate}
\end{lemma}

For $1<s<2$ and $h\ge h_0$, there exists a unique
level curve $H_*=h$ in $\rho\ge\rho_0$, given by
$\rho=\bar\rho(h,\theta)$.
Define $v(h,\theta)$ as the positive angular speed
along this curve and $\Omega(h)$ as its mean over one turn:
\begin{equation}\label{mean:definition}
 v(h,\theta):=
 \frac{\partial_\rho H_*(\bar\rho(h,\theta),\theta)}{\bar\rho(h,\theta)}
 \qquad \text{and} \qquad
 \Omega(h):=
 \left(\frac1{2\pi}\int_0^{2\pi}\frac{du}{v(h,u)}\right)^{-1}.
\end{equation}

\begin{lemma}[Mean angular speed]\label{lem:mean-velocity}
After increasing $h_0$ if necessary, for any $\theta \in \R$, $1<s<2$ and $h \ge h_0$, there hold
\begin{equation}\label{mean:bounds}
 v(h,\theta) \approx 1,\qquad \Omega(h)\approx 1,
\end{equation}
and
\begin{equation}\label{mean:derivative-bounds}
 \left|\partial_hv(h,\theta)\right|
                           \lesssim \bar\rho(h,0)^{-4}, \qquad  \Omega'(h)\gtrsim \bar\rho(h,0)^{-4}.
\end{equation}
Moreover, for any $\alpha,\beta\in\mathbb R$,
\begin{equation}\label{mean:primitive}
 \left|\int_\alpha^\beta\frac{du}{v(h,u)}
              -\frac{\beta-\alpha}{\Omega(h)}\right|\lesssim 1.
\end{equation}
\end{lemma}

\subsection{Variation of the collinear Hamiltonian}
We define $S_s(t)$ for $1\le s<2$ to express the
leading variation of $H_*$:
\begin{equation}\label{coefficient:matrix}
 S_s(t):=\sum_i a_i(t)\otimes a_i(t)
 -\frac{8\pi s^4}{3}
 \frac{a_1(t)^\perp\cdot\dot a_1(t)}{|a_1(t)|^2}
 \frac{a_1(t)\otimes a_1(t)}{|a_1(t)|^2}
\end{equation}
For $1<s<2$, we also define the coefficients $C_s$ and $C_*$ by
\begin{equation}\label{src:coefficient}
 C_s=\int_{-T_s/2}^{T_s/2}
 \left([S_s(t)]_{11}-[S_s(t)]_{22}\right)\,dt,\,\,\,\,\,C_*=\int_{\mathbb R}
 \bigl([S_1(t)]_{11}-[S_1(t)]_{22}\bigr)\,dt.
\end{equation}
\begin{lemma}[Hamiltonian variation]\label{lem:energy-change}\label{lem:coefficient-limit}
There exist universal constants $c>0$ and $\rho_0>0$ such that for any $1<s<2$, the following hold.
\begin{enumerate}
\item[(i)] For any $|t|\le T_s/2$, we have $\|S_s(t)\|\lesssim e^{-c|t|}$.

\item[(ii)] Let $X(t)$ be any solution of \eqref{ext:particle}.
For $|t|\le T_s/2$ with $\rho(t)\ge\rho_0$, we have
\begin{equation}\label{ext:angular-equation-error}
\begin{aligned}
 \left|\dot\theta(t)+
 \frac{\partial_\rho H_*(\rho(t),\theta(t))}{\rho(t)}
 \right|
 &\lesssim e^{-c|t|},\\
 \left|\frac d{dt}H_*(Z(t))\right|
 &\lesssim e^{-c|t|}\rho(t)^{-2}.
\end{aligned}
\end{equation}
Moreover,
\begin{equation}\label{ext:energy-leading}
 \frac d{dt}H_*(Z(t))
 =\frac3{4\pi^2s^2}
   \frac{X(t)^\perp\cdot S_s(t)X(t)}{|X(t)|^4}
 +O(e^{-c|t|}\rho(t)^{-3}).
\end{equation}

\item[(iii)] We have $\lim_{s \to 1^+}C_s=C_*>0$ and
$$ \int_{-T_s/2}^{T_s/2}[S_s(t)]_{12}\,dt=0.$$
\end{enumerate}
\end{lemma}
\begin{proof}
\noindent\textbf{Step 1: Matrix decay and the Hamiltonian difference.}
For any $\xi\in\mathbb R^2$, the identity
$(a\otimes a)\xi=a(a\cdot\xi)$ and \eqref{src:frame} give
\[
 \begin{aligned}
 \mathcal R(t)(a_i(t)\otimes a_i(t))\mathcal R(t)^T\xi=\mathcal R(t)a_i(t)
   \bigl(a_i(t)\cdot\mathcal R(t)^T\xi\bigr)
 =(b_i(t)\otimes b_i(t))\xi.
 \end{aligned}
\]
Thus,
$
 \mathcal R(t)(a_i(t)\otimes a_i(t))\mathcal R(t)^T
 =b_i(t)\otimes b_i(t)
$.
Since $\mathcal R(t)a_1(t)=|a_1(t)|e_1$ by
\eqref{src:frame}, we also have
\[
 \begin{aligned}
 \mathcal R(t)
 \frac{a_1(t)\otimes a_1(t)}{|a_1(t)|^2}
 \mathcal R(t)^T
 =\frac{(|a_1(t)|e_1)\otimes(|a_1(t)|e_1)}
         {|a_1(t)|^2}=e_1\otimes e_1.
 \end{aligned}
\]
Substitution in \eqref{coefficient:matrix} yields
\[
 \mathcal R(t)S_s(t)\mathcal R(t)^T
 =\sum_i b_i(t)\otimes b_i(t)
 -\frac{8\pi s^4}{3}
 \frac{a_1(t)^\perp\cdot\dot a_1(t)}{|a_1(t)|^2}
 (e_1\otimes e_1).
\]
For the collinear configuration $(se_1,0,-se_1)$,
\[
 (se_1)\otimes(se_1)+0 \otimes 0+(-se_1)\otimes(-se_1)
 =2s^2(e_1\otimes e_1).
\]
Using $\frac{8\pi s^4}{3}\frac3{4\pi s^2}=2s^2$,
we therefore rewrite the preceding identity as
\begin{equation}\label{coefficient:rotated-matrix}
 \begin{aligned}
 \mathcal R(t)S_s(t)\mathcal R(t)^T
 =&\sum_i b_i(t)\otimes b_i(t)-2s^2(e_1\otimes e_1)\\
 &-\frac{8\pi s^4}{3}
 \left(
 \frac{a_1(t)^\perp\cdot\dot a_1(t)}{|a_1(t)|^2}
 -\frac3{4\pi s^2}
 \right)(e_1\otimes e_1).
 \end{aligned}
\end{equation}
By \eqref{src:source-decay} and \eqref{src:speed-decay}, 
\begin{equation}\label{coefficient:comparison}
 \left\|\sum_i b_i(t)\otimes b_i(t)-2s^2(e_1\otimes e_1)\right\|
 +\left|
 \frac{a_1(t)^\perp\cdot\dot a_1(t)}{|a_1(t)|^2}
 -\frac3{4\pi s^2}
 \right|
 \lesssim e^{-c|t|}.
\end{equation}
Together with the orthogonality of $\mathcal R(t)$,
this gives \textup{(i)}.

For sufficiently large $\rho$ and $Z=(\rho\cos\theta,\rho\sin\theta)$,
Lemma~\ref{lem:logarithm}, \eqref{src:invariants}, and
\eqref{src:frame} give
\[
 \begin{aligned}
 \sum_i\log|Z-b_i(t)|
 =3\log\rho+\frac{s^2}{\rho^2}
 -\frac{\sum_i(Z\cdot b_i(t))^2}{\rho^4}
 +\sum_i\mathcal E(\rho,\theta,b_i(t)),
 \end{aligned}
\]
where we have used the following equalities from \eqref{src:invariants} and \eqref{src:frame}:
\[
 \sum_i b_i(t)=0,\qquad
 \sum_i|b_i(t)|^2=2s^2.
\]
Applying the same expansion to the collinear configuration $se_1,0,-se_1$, we obtain
\[
 \begin{aligned}
 &\log|Z-se_1|+\log|Z|+\log|Z+se_1|\\
 &=3\log\rho+\frac{s^2}{\rho^2}
 -\frac{2s^2Z_1^2}{\rho^4}+\mathcal E(\rho,\theta,se_1)
 +\mathcal E(\rho,\theta,0)
 +\mathcal E(\rho,\theta,-se_1).
 \end{aligned}
\]
Subtracting the two expansions and using
\eqref{ext:actual-hamiltonian} and \eqref{ext:reference} yields
\begin{equation}\label{ext:hamiltonian-difference}
 \begin{aligned}
 H(t,\rho,\theta)-H_*(\rho,\theta)
 ={}&\frac{\rho^2}{2}
 \left(
 \frac{a_1(t)^\perp\cdot\dot a_1(t)}{|a_1(t)|^2}
 -\frac3{4\pi s^2}
 \right)+\frac{\sum_i(Z\cdot b_i(t))^2-2s^2Z_1^2}
              {2\pi\rho^4}
 \\&+E(t,\rho,\theta),
 \end{aligned}
\end{equation}
where 
\[
 \begin{aligned}
 E(t,\rho,\theta):={}&-\frac1{2\pi}
 \Bigl[\sum_i\mathcal E(\rho,\theta,b_i(t))-\mathcal E(\rho,\theta,se_1)
 -\mathcal E(\rho,\theta,0)
 -\mathcal E(\rho,\theta,-se_1)\Bigr].
 \end{aligned}
\]
For $|t|\le T_s/2$, \eqref{src:source-decay} and
\eqref{ext:log-remainder} give
\begin{equation}\label{ext:hamiltonian-remainder}
 |E(t,\rho,\theta)|
 +\rho|\partial_\rho E(t,\rho,\theta)|
 +|\partial_\theta E(t,\rho,\theta)|
 \lesssim e^{-c|t|}\rho^{-3}.
\end{equation}

\medskip
\noindent\textbf{Step 2: The particle estimates.}
Differentiating \eqref{ext:hamiltonian-difference} yields
\begin{equation}\label{ext:radial-derivative-difference}
 \begin{aligned}
 \frac{\partial_\rho H-\partial_\rho H_*}{\rho}
 =&\frac{a_1(t)^\perp\cdot\dot a_1(t)}{|a_1(t)|^2}
 -\frac3{4\pi s^2}-\frac{\sum_i(Z\cdot b_i(t))^2-2s^2Z_1^2}{\pi\rho^6}
 \\&+\frac{\partial_\rho E(t,\rho,\theta)}{\rho}.
 \end{aligned}
\end{equation}
Using \eqref{coefficient:comparison}, we see that
\[
 \begin{aligned}
 \left|\sum_i(Z\cdot b_i(t))^2-2s^2Z_1^2\right|
 =\left|Z\cdot
 \left(\sum_i b_i(t)\otimes b_i(t)-2s^2(e_1\otimes e_1)\right)Z\right|\lesssim e^{-c|t|}\rho^2.
 \end{aligned}
\]
Combining \eqref{ext:radial-derivative-difference} with
\eqref{ext:polar-ode}, \eqref{ext:hamiltonian-remainder}
and \eqref{src:speed-decay}, we obtain
\[
 \begin{aligned}
 \left|\dot\theta(t)
 +\frac{\partial_\rho H_*(\rho(t),\theta(t))}{\rho(t)}\right|
 \lesssim e^{-c|t|},
 \end{aligned}
\]
which gives the first estimate in
\eqref{ext:angular-equation-error}.

To derive the second estimate in
\eqref{ext:angular-equation-error}, we differentiate
\eqref{ext:hamiltonian-difference} in $\theta$ and obtain
\begin{equation}\label{ext:angular-derivative-difference}
 \frac{\partial_\theta H-\partial_\theta H_*}{\rho}
 =\frac{
 Z^\perp\cdot
 \left(\sum_i b_i(t)\otimes b_i(t)-2s^2(e_1\otimes e_1)\right)Z}
 {\pi\rho^5}
 +\frac{\partial_\theta E(t,\rho,\theta)}{\rho}.
\end{equation}
Substituting \eqref{ext:radial-derivative-difference} and
\eqref{ext:angular-derivative-difference} into
\eqref{ext:energy-identity} yields
\[
 \begin{aligned}
 \frac d{dt}H_*(Z(t))
 ={}&\partial_\rho H_*
       \frac{\partial_\theta H-\partial_\theta H_*}{\rho}
 -\partial_\theta H_*
       \frac{\partial_\rho H-\partial_\rho H_*}{\rho}\\
 ={}&\frac{\partial_\rho H_*}{\pi\rho^5}
 Z^\perp\cdot
 \left(\sum_i b_i(t)\otimes b_i(t)-2s^2(e_1\otimes e_1)\right)Z\\
 &-\partial_\theta H_*
 \left(
 \frac{a_1(t)^\perp\cdot\dot a_1(t)}{|a_1(t)|^2}
 -\frac3{4\pi s^2}
 \right)+E_1,
 \end{aligned}
\]
where the new error term is 
$$
\begin{aligned}
E_1:=\frac{\partial_\rho H_*\,\partial_\theta E
              -\partial_\theta H_*\,\partial_\rho E}{\rho}+\frac{\partial_\theta H_*}{\pi\rho^6}
 \left(\sum_i(Z\cdot b_i(t))^2-2s^2Z_1^2\right).
\end{aligned}
$$
By \eqref{ext:reference-bounds},
\eqref{coefficient:comparison} and
\eqref{ext:hamiltonian-remainder}, we have $|E_1| \lesssim  e^{-c|t|}\rho^{-3}$ and
\begin{equation}\label{eq t derivetive of H * Z(t)}
 \begin{aligned}
 \frac d{dt}H_*(Z(t))
 ={}&\frac3{4\pi^2s^2\rho^4}
 Z^\perp\cdot
 \left(\sum_i b_i(t)\otimes b_i(t)-2s^2(e_1\otimes e_1)\right)Z\\
 &+\frac{s^2\sin2\theta}{\pi\rho^2}
 \left(
 \frac{a_1(t)^\perp\cdot\dot a_1(t)}{|a_1(t)|^2}
 -\frac3{4\pi s^2}
 \right)
 +O(e^{-c|t|}\rho^{-3}).
 \end{aligned}
\end{equation}
Since $Z^\perp\cdot(e_1\otimes e_1)Z
 =-\frac{\rho^2}{2}\sin2\theta$,
we have 
\[
\begin{aligned}
\frac{s^2\sin2\theta}{\pi\rho^2}
=-\frac{2s^2}{\pi\rho^4}
  Z^\perp\cdot(e_1\otimes e_1)Z=\frac3{4\pi^2s^2\rho^4}
  \left(-\frac{8\pi s^4}{3}\right)
  Z^\perp\cdot(e_1\otimes e_1)Z,
\end{aligned}
\]
which together with \eqref{eq t derivetive of H * Z(t)} and \eqref{coefficient:rotated-matrix} yield
\[
 \frac d{dt}H_*(Z(t))
 =\frac3{4\pi^2s^2\rho(t)^4}
 Z(t)^\perp\cdot\mathcal R(t)S_s(t)\mathcal R(t)^T Z(t)
 +O(e^{-c|t|}\rho(t)^{-3}).
\]
It follows from \eqref{ext:polar} that
\[
 Z(t)^\perp\cdot\mathcal R(t)S_s(t)\mathcal R(t)^T Z(t)
 =X(t)^\perp\cdot S_s(t)X(t),
 \qquad |X(t)|=\rho(t).
\]
Thus \eqref{ext:energy-leading} holds. Using \textup{(i)}
then gives the second estimate in
\eqref{ext:angular-equation-error}.

\medskip
\noindent\textbf{Step 3: The limiting coefficient.}
The initial data in Lemma~\ref{lem:vortex-invariants} depend
continuously on $s$ by
Lemma~\ref{lem:shape-function-properties}\textup{(ii)}.
Continuous dependence on initial data for ODE \eqref{src:cartesian},
the bound $|a_1(t)|=2r(t)\ge s>1$ in
Proposition~\ref{prop:source-family}\textup{(i)},
and \eqref{coefficient:matrix} give
\[
 S_s(t)\longrightarrow S_1(t)\qquad(s\to1^+)
 \quad\text{uniformly on every compact time interval}.
\]
Since $T_s\to\infty$ by Proposition~\ref{prop:source-family}\textup{(ii)},
the bound in \textup{(i)} implies
\[
 \|S_1(t)\|\lesssim e^{-c|t|},\qquad
 \|S_s(t)\|\mathbf1_{[-T_s/2,T_s/2]}(t)
 \lesssim e^{-c|t|}.
\]
The dominated convergence theorem applied to
\eqref{src:coefficient} gives $C_s\to C_*$.
The positivity $C_*>0$ follows from Lemma~\ref{lem:positive}.

By \eqref{eq:kernel} and Lemma~\ref{lem:vortex-invariants}, the three trajectories
$$(a_{1,1}(-t),-a_{1,2}(-t)), \quad
(a_{3,1}(-t),-a_{3,2}(-t)) \quad \text{and} \quad
(a_{2,1}(-t),-a_{2,2}(-t))$$
satisfy \eqref{src:cartesian} with the same initial data as
$a_1(t)$, $a_2(t)$ and $a_3(t)$.
Uniqueness of ODE system therefore gives
$a_1(-t)=(a_{1,1}(t),-a_{1,2}(t))$ and
$a_2(-t)=(a_{3,1}(t),-a_{3,2}(t))$.
Substitution in \eqref{coefficient:matrix} yields
$[S_s(-t)]_{12}=-[S_s(t)]_{12}$,
which proves \textup{(iii)}.
\end{proof}
\section{Hamiltonian change over one period}\label{sec:one-period}
We use the estimates of Section~\ref{sec:coordinates}
to compute the changes in $H_*$ and the particle's polar
angle during one period.

By Lemma~\ref{lem:energy-change}\textup{(iii)}, choose
$1<s<2$ sufficiently close to $1$ that $C_s\ge C_*/2$,
and set
\begin{equation}\label{cycle:parameters}
 L:=T_s^{1/8}.
\end{equation}
Recall that $L$ separates the transition interval $[-L,L]$
from the long waiting interval $[L,T_s-L]$ near the
collinear motion; see Figure~\ref{fig:intro-source-cycle}.
Then $L\to\infty$ as $s\to1^+$ by
Proposition~\ref{prop:source-family}\textup{(ii)}.

By \eqref{src:periodicity} and \eqref{src:periodicity-two},
the particle Hamiltonian has period $T_s$:
\begin{equation}\label{cycle:periodicity}
 H(t+T_s,Z)=H(t,Z).
\end{equation}
The estimates \eqref{ext:angular-equation-error} therefore
extend to all $t$ with $\rho(t)\ge\rho_0$, with
$e^{-c|t|}$ replaced by
$e^{-c\dist(t,T_s\mathbb Z)}$.
\begin{lemma}[One-period estimates]\label{lem:cycle}
For sufficiently large $L$ and any initial particle position satisfying $|X(-L)|\ge L/2$,
the solution of \eqref{ext:particle} exists on $[-L,T_s-L]$
and the following hold.
\begin{enumerate}
\item[(i)] For any $-L\le t\le T_s-L$,
\begin{equation}\label{cycle:confinement}
 |\rho(t)-\rho(-L)|\lesssim L^{-3}.
\end{equation}

\item[(ii)] The Hamiltonian change over one period satisfies
\begin{equation}\label{cycle:kick}
 \begin{aligned}
 \qquad H_*(Z(T_s-L))-H_*(Z(-L))
 =-\frac{3C_s}{4\pi^2s^2}
       \frac{X_1(-L)X_2(-L)}{\rho(-L)^4}+O(L^{-3}).
 \end{aligned}
\end{equation}

\item[(iii)] The angular displacement over one period satisfies
\begin{equation}\label{cycle:wait-angle}
 \theta(T_s-L)-\theta(-L)
 =-T_s\Omega\bigl(H_*(Z(T_s-L))\bigr)+O(1).
\end{equation}
\end{enumerate}
\end{lemma}
\begin{proof}
\noindent\textbf{Step 1: Radius and Hamiltonian bounds.}
As long as $\rho(-L)/2<\rho(t)<2\rho(-L)$, we have
$\rho(t)>L/4$.
For sufficiently large $L$, \eqref{ext:angular-equation-error}
and \eqref{cycle:periodicity} give
\begin{equation}\label{cycle:energy-control-proof}
 \begin{aligned}
 |H_*(Z(t))-H_*(Z(-L))|
 &\lesssim L^{-2}
 \int_{-L}^t e^{-c\dist(u,T_s\mathbb Z)}\,du
 \lesssim L^{-2}.
 \end{aligned}
\end{equation}
Equation \eqref{ext:energy-radius} therefore yields
\[
 |\rho(t)-\rho(-L)|
 \lesssim L^{-1}L^{-2}+L^{-3}
 \lesssim L^{-3}<\rho(-L)/2.
\]
By Lemma~\ref{lem:vortex-invariants}\textup{(i)},
\[
 |X(t)-a_i(t)|\ge\rho(t)-|a_i(t)|
 >L/4-2\sqrt2>0.
\]
These bounds extend the solution of \eqref{ext:particle}
to $[-L,T_s-L]$ and give \eqref{cycle:confinement}.

\medskip
\noindent\textbf{Step 2: The Hamiltonian change over one period.}
To prove \eqref{cycle:kick}, we write
\begin{equation}\label{eq sec4 difference H}
 \begin{aligned}
 &H_*(Z(T_s-L))-H_*(Z(-L))\\
 &\quad=\int_{-L}^{L}\frac d{dt}H_*(Z(t))\,dt
 +H_*(Z(T_s-L))-H_*(Z(L)).
 \end{aligned}
\end{equation}
For $L\le t\le T_s-L$, \eqref{ext:angular-equation-error},
\eqref{cycle:periodicity} and \eqref{cycle:confinement} give
\begin{equation}\label{cycle:energy-tail-proof}
 \begin{aligned}
 |H_*(Z(t))-H_*(Z(T_s-L))|
 \le\int_t^{T_s-L}\left|\frac d{du}H_*(Z(u))\right|\,du\lesssim e^{-cL}L^{-2}.
 \end{aligned}
\end{equation}
Since $\rho(t)\ge L/4$, using \eqref{ext:energy-leading}, \eqref{eq sec4 difference H} and
\eqref{cycle:energy-tail-proof}, we obtain
\begin{equation}\label{cycle:energy-integral-proof}
 \begin{aligned}
 H_*(Z(T_s-L))-H_*(Z(-L))=\frac3{4\pi^2s^2}\int_{-L}^{L}
 \frac{X(t)^\perp\cdot S_s(t)X(t)}{\rho(t)^4}\,dt
 +O(L^{-3}).
 \end{aligned}
\end{equation}

We replace $X(t)$ by $X(-L)$ in
\eqref{cycle:energy-integral-proof}.
For $-L\le t\le L$, \eqref{ext:particle},
\eqref{src:distance-bounds} and \eqref{cycle:confinement} give
\[
 |\dot X(t)|\lesssim L^{-1},\qquad
 |X(t)-X(-L)|\le\int_{-L}^t|\dot X(u)|\,du\lesssim1.
\]
Using \eqref{cycle:confinement}, a direct calculation yields
\[
 \begin{aligned}
 |X(t)^\perp\cdot S_s(t)X(t)
   -X(-L)^\perp\cdot S_s(t)X(-L)|
 &\lesssim \rho(-L)\|S_s(t)\||X(t)-X(-L)|\\
 &\lesssim \rho(-L)\|S_s(t)\|
 \end{aligned}
\]
and
\[
 \left|\rho(t)^{-4}-\rho(-L)^{-4}\right|
 \lesssim \rho(-L)^{-5}|\rho(t)-\rho(-L)|
 \lesssim \rho(-L)^{-5}L^{-3},
\]
which, together with
Lemma~\ref{lem:energy-change}\textup{(i)}, yields
\[
 \begin{aligned}
 \left|
 \frac{X(t)^\perp\cdot S_s(t)X(t)}{\rho(t)^4}
 -\frac{X(-L)^\perp\cdot S_s(t)X(-L)}{\rho(-L)^4}
 \right|
 \lesssim e^{-c|t|}L^{-3}.
 \end{aligned}
\]
Integrating over $[-L,L]$ gives
\begin{equation}\label{cycle:frozen-integral-proof}
 \begin{aligned}
 \int_{-L}^{L}
 \frac{X(t)^\perp\cdot S_s(t)X(t)}{\rho(t)^4}\,dt=\frac{X(-L)^\perp\cdot
 \int_{-L}^{L}S_s(t)\,dt X(-L)}
 {\rho(-L)^4}+O(L^{-3}).
 \end{aligned}
\end{equation}
To use \eqref{src:coefficient}, we extend the integral
to $[-T_s/2,T_s/2]$.
Lemma~\ref{lem:energy-change}\textup{(i)} gives
\begin{equation}\label{eq sec4 extend to Ts/2}
 \left\|
 \int_{-T_s/2}^{T_s/2}S_s(t)\,dt
 -\int_{-L}^{L}S_s(t)\,dt
 \right\|
 \lesssim\int_{|t|\ge L}e^{-c|t|}\,dt
 \lesssim e^{-cL}.
\end{equation}
Since $S_s(t)$ is symmetric by \eqref{coefficient:matrix}, using
\eqref{src:coefficient} and
Lemma~\ref{lem:energy-change}\textup{(iii)} we get
\[
 \begin{aligned}
 &X(-L)^\perp\cdot
 \int_{-T_s/2}^{T_s/2}S_s(t)\,dt X(-L)\\
 &\quad=-X_1(-L)X_2(-L)
 \int_{-T_s/2}^{T_s/2}
 \bigl([S_s(t)]_{11}-[S_s(t)]_{22}\bigr)\,dt\\
 &\quad=-C_sX_1(-L)X_2(-L).
 \end{aligned}
\]
Therefore, \eqref{cycle:kick} follows directly from \eqref{eq sec4 extend to Ts/2}, \eqref{cycle:frozen-integral-proof} and 
\eqref{cycle:energy-integral-proof}.

\medskip
\noindent\textbf{Step 3: The change in polar angle over one period.}
It follows from \eqref{cycle:confinement} and
Lemma~\ref{lem:reference-curves}\textup{(i)} that
$H_*(Z(t))\gtrsim L^2$ for $-L\le t\le T_s-L$.
Using \eqref{mean:primitive}, we get
\begin{equation}\label{cycle:angle-identity-proof}
 \begin{aligned}
 \frac{\theta(T_s-L)-\theta(-L)}
 {\Omega(H_*(Z(T_s-L)))}&=\int_{\theta(-L)}^{\theta(T_s-L)}
 \frac{du}{v(H_*(Z(T_s-L)),u)}+O(1)\\
 &=\int_{-L}^{T_s-L}
 \frac{\dot\theta(t)}{v(H_*(Z(T_s-L)),\theta(t))}\,dt+O(1)\\
 &=-T_s+\int_{-L}^{T_s-L}
 \frac{\dot\theta(t)+v(H_*(Z(T_s-L)),\theta(t))}
 {v(H_*(Z(T_s-L)),\theta(t))}\,dt+O(1).
 \end{aligned}
\end{equation}
By \eqref{mean:bounds}, it suffices to bound the last
integral to obtain \eqref{cycle:wait-angle}.
\[
 \begin{aligned}
 \dot\theta(t)+v(H_*(Z(T_s-L)),\theta(t))&=\dot\theta(t)+v(H_*(Z(t)),\theta(t))\\
 &\quad+v(H_*(Z(T_s-L)),\theta(t))
              -v(H_*(Z(t)),\theta(t)).
 \end{aligned}
\]
Lemma~\ref{lem:reference-curves}\textup{(iii)}
gives $\rho(t)=\bar\rho(H_*(Z(t)),\theta(t))$.
Thus \eqref{mean:definition}, \eqref{ext:angular-equation-error}
and \eqref{cycle:periodicity} imply
\begin{equation}\label{cycle:actual-angle-error-proof}
 \begin{aligned}
 \int_{-L}^{T_s-L}
 |\dot\theta(t)+v(H_*(Z(t)),\theta(t))|\,dt\lesssim\int_{-L}^{T_s-L}
 e^{-c\dist(t,T_s\mathbb Z)}\,dt\lesssim1.
 \end{aligned}
\end{equation}

Next, we estimate the difference
$v(H_*(Z(T_s-L)),\theta(t))-v(H_*(Z(t)),\theta(t))$.
For every $u$ between $H_*(Z(t))$ and $H_*(Z(T_s-L))$,
we have $u\gtrsim L^2$.
Lemma~\ref{lem:reference-curves}\textup{(i), (iii)} and
\eqref{mean:derivative-bounds} therefore give
\[
 |\partial_h v(u,\theta(t))|
 \lesssim \bar\rho(u,0)^{-4}\lesssim L^{-4},
\]
which leads to
\[
 \begin{aligned}
 |v(H_*(Z(T_s-L)),\theta(t))-v(H_*(Z(t)),\theta(t))|\lesssim L^{-4}|H_*(Z(T_s-L))-H_*(Z(t))|.
 \end{aligned}
\]
Combining this with  \eqref{cycle:parameters}, \eqref{cycle:energy-control-proof} and \eqref{cycle:energy-tail-proof}, we obtain
\begin{equation}\label{cycle:velocity-error-integral-proof}
 \begin{aligned}
 &\int_{-L}^{T_s-L}
 |v(H_*(Z(t)),\theta(t))
   -v(H_*(Z(T_s-L)),\theta(t))|\,dt\\
 &\quad\lesssim L^{-4}\int_{-L}^{L}
 |H_*(Z(t))-H_*(Z(T_s-L))|\,dt\\
 &\qquad+L^{-4}\int_L^{T_s-L}
 |H_*(Z(t))-H_*(Z(T_s-L))|\,dt\\
 &\quad\lesssim L^{-5}+T_se^{-cL}L^{-6}
 \lesssim1.
 \end{aligned}
\end{equation}
Substituting \eqref{mean:bounds},
\eqref{cycle:actual-angle-error-proof} and
\eqref{cycle:velocity-error-integral-proof} into
\eqref{cycle:angle-identity-proof}, we obtain
\eqref{cycle:wait-angle}.
\end{proof}

\section{Iteration and finite growth}\label{sec:arc-selection}
We now iterate the one-period estimates to obtain a
particle reaching a prescribed large radius.
Rescaling the motion yields the following proposition.

\begin{proposition}\label{prop:passive}
For every $N>1$, there exist $T>0$, a solution
$a_1(t),a_2(t),a_3(t)$ of \eqref{src:cartesian}, and a solution
$X(t)$ of \eqref{ext:particle} on $[0,T]$ such that
$a_i(0),X(0)\in B_{3/4}$ for $i=1,2,3$.
The mutual distances of these four paths have a positive lower
bound on $[0,T]$ and
\[
 \min_{1\le i\le3}|X(T)-a_i(T)|>N+2.
\]
\end{proposition}

\subsection{Selection after one period}
Fix $N>1$. We take $s>1$ sufficiently close to $1$, depending on $N$, and define $L$ by \eqref{cycle:parameters}.
All implicit constants in this section may depend on $N$
but are independent of $s$.

For $|z|\ge L/2$, let $Z(t)$ be the particle trajectory
in rotating coordinates with $Z(-L)=z$.
By Lemma~\ref{lem:cycle}, we define the one-period map
\begin{equation}\label{select:return-map}
 \mathcal P(z):=Z(T_s-L).
\end{equation}
By \eqref{cycle:periodicity}, for $n=0,1,\ldots$,
$\mathcal P^n(z)=Z(nT_s-L)$ whenever the iterates are defined.
To increase $H_*$ under iteration, we choose initial
angles making the leading term in \eqref{cycle:kick}
positive; see Figure~\ref{fig:intro-sectors}.
More specifically, define $\alpha\in\mathbb R$ and
$U\subset\mathbb R$ by
\begin{equation}\label{select:angle-interval}
 \mathcal R(-L)^T(\cos\alpha,\sin\alpha)
 =\frac1{\sqrt2}(1,-1),\qquad
 U:=[\alpha,\alpha+\pi/6].
\end{equation}
For $z=(\rho\cos\theta,\rho\sin\theta)$, \eqref{ext:polar} and
\eqref{select:angle-interval} give
\[
 -\frac{X_1(-L)X_2(-L)}{\rho^4}
 =\frac{\cos(2(\theta-\alpha))}{2\rho^2}.
\]
Thus \eqref{cycle:kick} becomes
\begin{equation}\label{select:sector-gain}
 H_*(\mathcal P(z))-H_*(z)
 =\frac{3C_s}{8\pi^2s^2\rho^2}\cos(2(\theta-\alpha))
 +O(L^{-3}).
\end{equation}

To repeat the increase, we select a subcurve with
angle range $U$ and Hamiltonian variation at most
$L^{-3}$; see Figure~\ref{fig:intro-arc-evolution}.

\begin{lemma}[Return of a curve]\label{lem:arc-return}
For any $N>1$, there exist $\delta_N\in(0,1)$ and $c_N>0$
such that the following holds for any $1<s<1+\delta_N$:
Consider a curve
\[
 \Gamma(q)=(\rho(q)\cos\varphi(q),\rho(q)\sin\varphi(q)),\qquad 0\le q\le1,
\]
where $\rho,\varphi$ are continuous and satisfy
$L/2\le\rho(q)\le2(N+3)L$,
\begin{equation}\label{select:angle-range}
 \varphi(q)\in U,\qquad
 \varphi(0)=\alpha,\qquad\varphi(1)=\alpha+\pi/6,
\end{equation}
and
\begin{equation}\label{select:energy-width}
 \max_{p,q\in[0,1]}
 |H_*(\Gamma(p))-H_*(\Gamma(q))|\le L^{-3}.
\end{equation}
Then $\mathcal P(\Gamma([0,1]))$ contains a continuous curve that,
after reparametrization on $[0,1]$, satisfies
\eqref{select:angle-range} and \eqref{select:energy-width}.
Moreover,
\begin{equation}\label{select:positive-gain}
 H_*(\mathcal P(\Gamma(q)))-H_*(\Gamma(q))\ge c_N L^{-2},
 \qquad 0\le q\le1.
\end{equation}
\end{lemma}
\begin{proof}
Fix $N>1$. All constants below may depend on $N$,
but are independent of $s$ and $\Gamma$.

\noindent\textbf{Step 1: Hamiltonian increase.}
For $s>1$ sufficiently close to $1$,
Lemma~\ref{lem:energy-change}\textup{(iii)} gives $C_s\ge C_*/2$.
For $\varphi\in U$, we have $\cos(2(\varphi-\alpha))\ge1/2$.
Thus \eqref{select:sector-gain} and
$L/2\le\rho(q)\le2(N+3)L$ give \eqref{select:positive-gain}
with some $c_N>0$.
By \eqref{ext:energy-radius} and
\eqref{select:energy-width},
\begin{equation}\label{eq sec5 rho p- rho q}
 |\rho(p)-\rho(q)|
 \lesssim L^{-1}L^{-3}+L^{-3}\lesssim L^{-3},
 \qquad p,q\in[0,1].
\end{equation}
Using \eqref{select:sector-gain}, \eqref{eq sec5 rho p- rho q}, \eqref{select:angle-range} and
\eqref{select:energy-width}, we obtain
\begin{equation}\label{select:endpoint-energy-gap-proof}
 \begin{aligned}
 H_*(\mathcal P(\Gamma(0)))-H_*(\mathcal P(\Gamma(1)))=\frac{3C_s}{8\pi^2s^2}
 \left(\frac1{\rho(0)^2}-\frac1{2\rho(1)^2}\right)
 +O(L^{-3})\gtrsim L^{-2}.
 \end{aligned}
\end{equation}

\medskip
\noindent\textbf{Step 2: Change in polar angle.}
For each $q\in[0,1]$, consider the solution $X(t)$ of
\eqref{ext:particle} with
$X(-L)=\mathcal R(-L)^T\Gamma(q)$ and set
$Z(t)=\mathcal R(t)X(t)$ as in \eqref{ext:polar}.
Let $\theta(q,t)$ be the polar angle of $Z(t)$,
counting full turns, with $\theta(q,-L)=\varphi(q)$.
Continuous dependence on initial data implies
that $q\mapsto\theta(q,T_s-L)$ is continuous.
By \eqref{cycle:wait-angle},
\begin{equation}\label{select:return-angle-proof}
 \theta(q,T_s-L)=\varphi(q)
 -T_s\Omega(H_*(\mathcal P(\Gamma(q))))+O(1),
\end{equation}
which implies
\[
 \begin{aligned}
 \theta(1,T_s-L)-\theta(0,T_s-L)=\frac\pi6+T_s\bigl[
 \Omega(H_*(\mathcal P(\Gamma(0))))
 -\Omega(H_*(\mathcal P(\Gamma(1))))\bigr]+O(1).
 \end{aligned}
\]
Using  \eqref{cycle:confinement} and
Lemma~\ref{lem:reference-curves}\textup{(i)}, we obtain
$H_*(\mathcal P(\Gamma(q)))\approx L^2$.
For any $h\approx L^2$,
Lemma~\ref{lem:reference-curves}\textup{(i), (iii)} and
\eqref{mean:derivative-bounds} give
\begin{equation}\label{select:frequency-slope-proof}
 \bar\rho(h,0)\lesssim L,\qquad \Omega'(h)\gtrsim L^{-4}.
\end{equation}
Therefore, 
\eqref{select:endpoint-energy-gap-proof},
\eqref{select:frequency-slope-proof}, \eqref{cycle:parameters} and the mean value theorem give
\begin{equation}\label{eq angle growth}
 \theta(1,T_s-L)-\theta(0,T_s-L)
 \ge cT_sL^{-6}-C=cL^2-C>4\pi.
\end{equation}

\medskip
\noindent\textbf{Step 3: Selection of the subcurve.}
By continuity and \eqref{eq angle growth}, there exist $k\in\mathbb Z$ and $0\le q_0<q_1\le1$ such that
\[
 \begin{gathered}
 \theta(q_0,T_s-L)=\alpha+2\pi k,\qquad
 \theta(q_1,T_s-L)=\alpha+\pi/6+2\pi k,\\
 \theta(q,T_s-L)\in U+2\pi k,\qquad q\in[q_0,q_1].
 \end{gathered}
\]
For $0\le u\le1$, define
\[
 \begin{aligned}
 \widetilde\Gamma(u):
 &=\mathcal P\bigl(\Gamma(q_0+(q_1-q_0)u)\bigr),\\
 \widetilde\varphi(u):
 &=\theta(q_0+(q_1-q_0)u,T_s-L)-2\pi k.
 \end{aligned}
\]
Then $\widetilde\Gamma$ satisfies
\eqref{select:angle-range}.
By  \eqref{select:angle-range}, \eqref{select:return-angle-proof}  and $\widetilde\varphi([0,1])\subset U$, we have
\[
 T_s|\Omega(H_*(\widetilde\Gamma(u)))
       -\Omega(H_*(\widetilde\Gamma(v)))|\le C,
 \qquad u,v\in[0,1].
\]
The mean value theorem and \eqref{select:frequency-slope-proof} imply that
\[
 \max_{u,v\in[0,1]}
 |H_*(\widetilde\Gamma(u))-H_*(\widetilde\Gamma(v))|
 \le\frac{CL^4}{T_s}\le L^{-3}
\]
for sufficiently large $L$, so that $\widetilde\Gamma$
satisfies \eqref{select:energy-width}.

Since $L\to\infty$ and $C_s\to C_*>0$ as $s\to1^+$,
choosing $\delta_N\in(0,1)$ sufficiently small ensures all the
required estimates for every $1<s<1+\delta_N$.
\end{proof}

\subsection{Iteration and rescaling}
\begin{proof}[Proof of Proposition~\ref{prop:passive}]
Fix $N>1$ and choose $1<s<1+\delta_N$ sufficiently close to $1$,
where $\delta_N$ is given by Lemma~\ref{lem:arc-return}.
For $0\le q\le1$, define
\[
 \varphi(q):=\alpha+\frac{\pi q}{6},\qquad
 \rho(q):=\bar\rho(H_*(L,0),\varphi(q)).
\]
The initial curve is
\[
 \Gamma_0=
 \left\{(\rho(q)\cos\varphi(q),\rho(q)\sin\varphi(q)):
 0\le q\le1\right\}.
\]
For any $z\in\Gamma_0$, write
$z=(|z|\cos\theta,|z|\sin\theta)$ with $\theta\in U$.
Lemma~\ref{lem:reference-curves}\textup{(i), (iii)} gives
$H_*(z)=H_*(L,0)$ and $|z|\approx L$.
Using \eqref{ext:energy-radius}, we obtain
\begin{equation}\label{select:initial-radius-proof}
 |z|=\bar\rho(H_*(L,0),\theta)=L+O(L^{-3}),
 \qquad \theta\in U.
\end{equation}
Thus, $\Gamma_0$ with the parametrization above satisfies
the assumptions of Lemma~\ref{lem:arc-return}.

Apply Lemma~\ref{lem:arc-return} repeatedly, selecting
$\Gamma_{n+1}\subset\mathcal P(\Gamma_n)$ as long as
$\max_{z\in\Gamma_n}|z|<2(N+3)L$.
We verify inductively that every selected curve satisfies
\begin{equation}\label{eq sec5 induction}
 |z|\ge L/2,\qquad H_*(z)\ge H_*(L,0),
 \qquad z\in\Gamma_n.
\end{equation}
These inequalities hold on $\Gamma_0$ by its construction.
Suppose they hold on $\Gamma_n$.
For any $y\in\Gamma_{n+1}$, there exists $z\in\Gamma_n$
such that $y=\mathcal P(z)$.
By \eqref{cycle:confinement},
\begin{equation}\label{eq sec5 low bound L/4}
 |y|\ge |z|-CL^{-3}\ge L/2-CL^{-3}\ge L/4.
\end{equation}
Write $y=(|y|\cos\theta,|y|\sin\theta)$ with $\theta\in U$.
Using \eqref{mean:level-radius} and \eqref{select:positive-gain}, we obtain
\[
 H_*(|y|,\theta)\ge H_*(L,0)
 =H_*(\bar\rho(H_*(L,0),\theta),\theta).
\]
Using \eqref{ext:reference-bounds}, we see that
$\partial_\rho H_*(\rho,\theta)>0$ for $\rho\ge L/4$.
Thus, \eqref{eq sec5 low bound L/4} and
\eqref{select:initial-radius-proof} yield
\[
 |y|\ge\bar\rho(H_*(L,0),\theta)
 \ge L-CL^{-3}>L/2.
\]
Thus \eqref{eq sec5 induction} holds on $\Gamma_{n+1}$.

If the iteration never stops, then
$L/2\le |z|<2(N+3)L$ on every $\Gamma_n$.
By \eqref{select:positive-gain} and
Lemma~\ref{lem:reference-curves}\textup{(i)}, we obtain
\[
 H_*(L,0)+nc_NL^{-2}\le H_*(z)\le C_NL^2,
 \qquad z\in\Gamma_n,\quad n\ge0,
\]
which is impossible as $n\to\infty$.
Hence, for some finite $m>0$, $\Gamma_m$ contains a point
of radius at least $2(N+3)L$.
Taking successive preimages of this point in the selected curves,
we obtain a particle trajectory on $[-L,mT_s-L]$ satisfying
\[
 |X(mT_s-L)|\ge2(N+3)L,\qquad
 |X(t)|\ge L/4,\quad -L\le t\le mT_s-L
\]
with $Z(-L)\in\Gamma_0$. By \eqref{src:distance-bounds}, we have
\[
 |X(t)-a_i(t)|\ge L/4-2\sqrt2>0,
 \qquad -L\le t\le mT_s-L,\quad i=1,2,3.
\]
Together with the vortex separation in \eqref{src:distance-bounds},
this gives a positive lower bound for all mutual distances.

Set $\lambda=(2L)^{-1}$, $T=\lambda^2mT_s$, and
\[
 \widetilde a_i(t)=\lambda a_i(t/\lambda^2-L),\qquad
 \widetilde X(t)=\lambda X(t/\lambda^2-L),
 \qquad 0\le t\le T.
\]
Since $K(\lambda z)=\lambda^{-1}K(z)$, these paths satisfy
\eqref{src:cartesian} and \eqref{ext:particle}.
By \eqref{src:distance-bounds} and \eqref{select:initial-radius-proof},
\[
 \max_i|\widetilde a_i(0)|\le\frac{\sqrt2}{L}<\frac34,
 \qquad
 |\widetilde X(0)|=\frac12+O(L^{-4})<\frac34,
\]
and
\[
 \min_i|\widetilde X(T)-\widetilde a_i(T)|
 \ge N+3-\frac{\sqrt2}{L}>N+2.
\]
This scaling preserves the positive lower bound
for the mutual distances.
\end{proof}

\section{Proof of the main theorem}\label{sec:patch}
We give the passive particle a small positive circulation
and use point-vortex approximation to obtain a smooth
simply connected patch with the required diameter.
\begin{proof}[Proof of Theorem~\ref{thm:main}]
Fix $N>1$ and take the initial positions $a_i(0)$, $1\le i\le3$,
$X(0)$ and the time $T$ from Proposition~\ref{prop:passive}.
Consider the four-vortex system
\[
 \dot a_i(t)=\sum_{\substack{1\le j\le4\\j\ne i}}
 \gamma_jK(a_i(t)-a_j(t)),\qquad i=1,2,3,4,
\]
where $\gamma_1=\gamma_2=\gamma_3=1$ and $\gamma_4=\gamma>0$,
with the first three initial positions unchanged and $a_4(0)=X(0)$.
At $\gamma=0$, this is the system in Proposition~\ref{prop:passive}.
For sufficiently small $\gamma\in(0,1)$, continuous
dependence on $\gamma$ gives four separated paths
on $[0,T]$ satisfying
\begin{equation}\label{patch:point-distance}
 \min_{1\le i\le3}|a_4(T)-a_i(T)|>N+1.
\end{equation}

For sufficiently small $\varepsilon>0$, the four disks in
\[
 E_0:=\bigcup_{i=1}^4
 B_{\varepsilon\sqrt{\gamma_i}}(a_i(0))
\]
are disjoint and contained in $B_{3/4}$.
Joining them in a tree by three nonintersecting arcs and taking
a sufficiently thin smooth neighborhood, we obtain a simply
connected domain $D_0$ with $C^\infty$ boundary such that
\begin{equation}\label{patch:corridor-area}
 \overline{E_0}\subset D_0\subset B_1,\qquad
 |D_0\setminus E_0|<\varepsilon^{100}.
\end{equation}
Define the initial vorticity by
\[
 \begin{aligned}
 \omega_0(x)
 &:=\frac1{\pi\varepsilon^2}\mathbf1_{D_0}(x)\\
 &=\frac1{\pi\varepsilon^2}\sum_{i=1}^4
 \mathbf1_{B_1}\left(\frac{x-a_i(0)}
                            {\varepsilon\sqrt{\gamma_i}}\right)
 +\frac1{\pi\varepsilon^2}\mathbf1_{D_0\setminus E_0}(x).
 \end{aligned}
\]
The four disk terms have circulations $\gamma_i$, and
\eqref{patch:corridor-area} gives
\[
 \left\|\frac1{\pi\varepsilon^2}
 \mathbf1_{D_0\setminus E_0}\right\|_{L^1}
 \lesssim\varepsilon^{98}.
\]
Let $\omega(t)$ be the Euler solution with initial vorticity
$\omega_0$. Yudovich transport
\cite[Theorems 1.6 and 3.3]{CrippaStefani} gives
\[
 \omega(t)=\frac1{\pi\varepsilon^2}\mathbf1_{D_t},
\]
where $D_t$ is the image of $D_0$ under the Euler flow.
For sufficiently small $\varepsilon$, the point-vortex
approximation \cite[Theorem 1.1]{Guo} gives
\[
 D_t\cap B_{1/4}(a_i(t))\ne\varnothing,
 \qquad i=1,2,3,4,\quad 0\le t\le T.
\]
Consequently, \eqref{patch:point-distance} gives
\begin{equation}\label{patch:diameter}
 \operatorname{diam}D_T>N+1-\frac12>N.
\end{equation}

The rescaled solution $\pi\varepsilon^2\omega(\pi\varepsilon^2t,x)$
has initial vorticity $\mathbf1_{D_0}$ and equals
$\mathbf1_{D_T}$ at $T_N=\frac{T}{\pi\varepsilon^2}$.
The Euler flow preserves simple connectivity
\cite[Section 1.3]{MajdaBertozzi}. Boundary regularity persists
at every finite order by \cite[Corollary B]{Chemin}.
\end{proof}

\appendix

\section{Properties of \texorpdfstring{$F_s$}{F_s}}\label{app:shape-function}
\begin{lemma}[Properties of $F_s$]\label{lem:shape-function-properties}
Let $F_s$ be defined by \eqref{src:shape-function}, then
\begin{equation}\label{src:shape-function-alt}
 F_s(x,y)=(x^2+y^2)\bigl[(3s^2-4x^2-4y^2)^2
                  +12(s^2-x^2-y^2)y^2\bigr],
\end{equation}
\begin{equation}\label{src:shape-derivatives}
 \partial_xF_s(x,y)
 =6x\bigl[(3s^2-4x^2-4y^2)(s^2-4x^2-4y^2)
              +4(s^2-2x^2-2y^2)y^2\bigr]
\end{equation}
and
\begin{equation}\label{src:shape-derivative-y}
 \partial_yF_s(x,y)
 =6y(s^2-2x^2-2y^2)(3s^2-6x^2-2y^2).
\end{equation}
The following properties hold.
\begin{enumerate}
\item[(i)] For any $s\ge1$ and $(x,y)\in\mathbb R^2$,
$F_s(-x,y)=F_s(x,y)=F_s(x,-y)$.

\item[(ii)] For any $s\ge1$ and $|y|<s/2$, the equation
\[
 F_s(0,y)=y^2(3s^2-2y^2)^2=1
\]
has exactly two solutions, denoted by $y_s(0)$ and $-y_s(0)$,
where $y_s(0)$ depends smoothly on $s$ and satisfies
$0<y_s(0)<s/2$.

\item[(iii)] For every $s>1$, the level curve
\begin{equation}\label{def level curve}
 \{(x,y)\in\mathbb R^2:F_s(x,y)=1,\ x^2+y^2<s^2/4\}
\end{equation}
is a connected smooth closed curve on which $|\nabla F_s(x,y)|>0$.

\item[(iv)] Let $1<s\le 2$, then on the level curve
\eqref{def level curve}, there hold
\begin{equation}\label{src:radial-bounds}
 \frac1{9s^4}\le x^2+y^2<\frac{s^2}{4}
\end{equation}
and
\begin{equation}\label{fs:level-displacement}
 s/2-|x|\approx |y|+\sqrt{s-1}.
\end{equation}
The derivatives satisfy
\begin{equation}\label{fs:level-derivative-x}
 x\partial_xF_s(x,y)
 \approx x^2\bigl(|y|+\sqrt{s-1}\bigr)
\end{equation}
and
\begin{equation}\label{fs:level-derivative-y}
 y\partial_yF_s(x,y)\approx y^2.
\end{equation}
\end{enumerate}
All implicit constants above are independent of $s,x,y$ in the
stated ranges.
\end{lemma}
\begin{proof}
\noindent\textbf{Step 1: Algebraic identities and the initial height.}
The identities \eqref{src:shape-function-alt}, \eqref{src:shape-derivatives},
\eqref{src:shape-derivative-y} and property (i) follow directly from
\eqref{src:shape-function}.

For any $s\ge1$, we see that $F_s(0,0)=0$, $F_s(0,s/2)=\frac{25s^6}{16}>1$
and
\[
 \partial_yF_s(0,y)
 =6y(3s^2-2y^2)(s^2-2y^2)>0,
 \qquad 0<y<s/2,
\]
so there exists a unique $0<y_s(0)<s/2$ such that $F_s(0,y_s(0))=1$.
The case $y<0$ follows from (i).
Since $\partial_yF_s(0,y_s(0))>0$, the implicit function theorem
applied to $(s,y)\mapsto F_s(0,y)-1$ gives smooth dependence of
$y_s(0)$ on $s$, which proves (ii).

\medskip
\noindent\textbf{Step 2: The level curve.}
Let $s>1$ and $x^2+y^2<s^2/4$.
Then
\begin{equation}\label{eq s1}
 s^2-4x^2-4y^2>0,\quad
 3s^2-4x^2-4y^2>2s^2,\quad s^2-2x^2-2y^2>\frac{s^2}{2}
\end{equation}
and
\begin{equation}\label{eq s2}
 3s^2-6x^2-2y^2
 \ge3s^2-6(x^2+y^2)>3s^2/2.
\end{equation}
Equations \eqref{src:shape-derivatives} and
\eqref{src:shape-derivative-y} therefore give
\[
 x\partial_xF_s(x,y)>0\quad\text{if }\,x\ne0,
 \qquad
 y\partial_yF_s(x,y)>0\quad\text{if }\,y\ne0.
\]
Consequently, on the region $0<x^2+y^2<s^2/4$,
\begin{equation}\label{eq app radial derivative}
 x\partial_xF_s(x,y)+y\partial_yF_s(x,y)>0,
 \qquad \text{and}\qquad \nabla F_s(x,y)\ne0.
\end{equation}
On the boundary $x^2+y^2=s^2/4$, \eqref{src:shape-function-alt}
gives $F_s(x,y)=s^6+\frac94s^4y^2>1$.

For each $\theta\in[0,2\pi]$, \eqref{eq app radial derivative}
shows that $\rho\mapsto F_s(\rho\cos\theta,\rho\sin\theta)$
is strictly increasing on $(0,s/2)$.
Since $F_s(0,0)=0$ and $F_s>1$ on the boundary $\rho=s/2$,
there exists a unique $0<\rho<s/2$ such that
$F_s(\rho\cos\theta,\rho\sin\theta)=1$.
Moreover, for any $s>1$, the smoothness of this level curve follows
directly from \eqref{eq app radial derivative} and the implicit
function theorem, which proves (iii).

\medskip
\noindent\textbf{Step 3: Uniform estimates on the level curve.}
Now it remains to prove (iv).
Let $1<s\le 2$ and consider $(x,y)$ on the level curve
\eqref{def level curve}.

By \eqref{src:shape-function},
\[
 1=F_s(x,y)
 \le(x^2+y^2)(3s^2-2x^2-2y^2)^2
 \le9s^4(x^2+y^2),
\]
which proves \eqref{src:radial-bounds}.
Since $1<s\le 2$, \eqref{eq s1} and \eqref{eq s2} yield
\[
 y\partial_yF_s(x,y)
 =6y^2(s^2-2x^2-2y^2)(3s^2-6x^2-2y^2)
 \approx y^2,
\]
which proves \eqref{fs:level-derivative-y}.
Similarly, \eqref{src:shape-derivatives} gives
\[
 \begin{split}
 x\partial_xF_s(x,y)
 &=6x^2\bigl[(3s^2-4x^2-4y^2)(s^2-4x^2-4y^2)
              +4(s^2-2x^2-2y^2)y^2\bigr]\\
 &\approx x^2\bigl[(s^2/4-x^2-y^2)+y^2\bigr]
 =x^2(s/2+|x|)(s/2-|x|)\\
 &\approx x^2(s/2-|x|).
 \end{split}
\]
Since $F_s(x,y)$ is even in each variable, it suffices to prove
\eqref{fs:level-displacement} for $x\ge0$ and $y\ge0$.
Using $F_s(x,y)=1$, we write
\begin{equation}\label{eq app split of F}
 F_s(s/2,0)-F_s(x,0)
 =F_s(s/2,0)-1+F_s(x,y)-F_s(x,0).
\end{equation}
On one hand, \eqref{src:shape-function-alt} gives
\[
 F_s(s/2,0)-F_s(x,0)
 =16(s^2-x^2)(s/2-x)^2(s/2+x)^2\approx(s/2-x)^2.
\]
On the other hand, \eqref{src:shape-derivative-y},
\eqref{eq s1} and \eqref{eq s2} yield
\[
 \begin{split}
 F_s(x,y)-F_s(x,0)
 &=\int_0^y6u(s^2-2x^2-2u^2)(3s^2-6x^2-2u^2)\,du\\
 &\approx\int_0^y u\,du=\frac{y^2}{2}.
 \end{split}
\]
Note that $F_s(s/2,0)-1=s^6-1\approx s-1$.
Using \eqref{eq app split of F}, we obtain
$(s/2-x)^2\approx s-1+y^2$, which proves \eqref{fs:level-displacement}
and hence \eqref{fs:level-derivative-x}, completing (iv).
\end{proof}

\section{Calculation and positivity of the limiting coefficient}
\label{app:coefficient}

\begin{lemma}[Positivity of the limiting coefficient]\label{lem:positive}
The coefficient $C_*$ defined in \eqref{src:coefficient} satisfies
\begin{equation}\label{coefficient:positive-bound}
 C_*>0.
\end{equation}
\end{lemma}

Before its proof, we first estimate the s-vortex motion with $s=1$.

\begin{lemma}\label{lem:coefficient-trajectory}
Let $\bar a_i(t)$ be the solution in
Lemma~\ref{lem:vortex-invariants} with $s=1$. Denote its coordinates
in \eqref{src:geometry} by $\bar r(t),\bar x(t),\bar y(t)$.
\begin{enumerate}
\item[(i)] For any $t\in\mathbb R$, we have
$\bar x(t)^2+\bar y(t)^2\le1/4$, $\dot{\bar x}(t)<0$ and
$\bar y(t)=g(\bar x(t))>0$, where
\begin{equation}\label{coefficient:g}
 g(x):=\sqrt{1-\frac52x^2
 -\frac{\sqrt3}{2}\sqrt{(1-x^2)(1-3x^2)}},
 \qquad |x|<\frac12.
\end{equation}
Moreover,
\begin{equation}\label{coefficient:coordinate-endpoints}
 \lim_{t\to-\infty}\bar x(t)=\frac12,
 \qquad
 \lim_{t\to+\infty}\bar x(t)=-\frac12.
\end{equation}

\item[(ii)] There is a smooth function $\Theta:(-1/2,1/2)\to\mathbb R$
with $\Theta(0)=0$ such that
\begin{equation}\label{coefficient:direction}
 \frac{\bar a_1(t)}{|\bar a_1(t)|}
 =\bigl(\cos\Theta(\bar x(t)),\sin\Theta(\bar x(t))\bigr).
\end{equation}
This angle is given by
\begin{equation}\label{coefficient:angle}
 \Theta(x)
 =-2\int_0^x
 \frac{\bigl(\xi^2+g(\xi)^2\bigr)
       \bigl(3-4\xi^2-2g(\xi)^2\bigr)}
 {g(\xi)\bigl(1-2\xi^2-2g(\xi)^2\bigr)
        \bigl(3-6\xi^2-2g(\xi)^2\bigr)}\,d\xi.
\end{equation}
\end{enumerate}
\end{lemma}
\begin{proof}
\noindent\textbf{Step 1: The coordinate relation.}
Proposition~\ref{prop:source-family}\textup{(i)},
Lemma~\ref{lem:shape-equations}\textup{(i)}, \eqref{src:frame} and continuous dependence
on $s$ give $\bar y(0)>0$, $F_1(\bar x(t),\bar y(t))=1$ and
\[
 \bar x(t)^2+\bar y(t)^2\le\frac14,
 \qquad t\in\mathbb R.
\]
Therefore, if $\bar y(t)=0$, then
\[
 0=F_1(\bar x(t),0)-1
 =\bigl(\bar x(t)^2-1\bigr)\bigl(4\bar x(t)^2-1\bigr)^2,
\]
which yields $\bar x(t)=\pm1/2$.
The identity $\nabla F_1(\pm1/2,0)=0$ implies that
$(\pm1/2,0)$ are stagnation points of the coordinate equations,
contradicting $\bar y(0)>0$. Hence $\bar y(t)>0$ for all
$t\in\mathbb R$.

We next show that, in the region $y>0$ and $x^2+y^2\le1/4$,
the equation $F_1(x,y)=1$ is equivalent to $y=g(x)$.
Indeed,
\[
 F_1(x,y)-1
 =(x^2+y^2-1)
 \bigl[4y^4+(20x^2-8)y^2+(1-4x^2)^2\bigr],
\]
which implies
\[
 4y^4+(20x^2-8)y^2+(1-4x^2)^2=0.
\]
Solving this quadratic equation for $y^2$ gives
\[
 y^2=1-\frac52x^2
 \pm\frac{\sqrt3}{2}\sqrt{(1-x^2)(1-3x^2)}.
\]
Since $1-\frac52x^2>\frac14$, the root with the plus sign exceeds
$1/4$, so we have
\[
 y^2=1-\frac52x^2
 -\frac{\sqrt3}{2}\sqrt{(1-x^2)(1-3x^2)},
\]
which proves $y=g(x)$ since $y>0$.

Since $\bar x(t)^2+\bar y(t)^2\le1/4$,
Lemma~\ref{lem:shape-equations}\textup{(iii)} and
\eqref{src:shape-derivative-y} give $\dot{\bar x}(t)<0$.
Boundedness and monotonicity of $\bar x(t)$, together with
$\bar y(t)=g(\bar x(t))$, imply that the limits of
$(\bar x(t),\bar y(t))$ exist as $t\to\pm\infty$,
and these limits must be stagnation points.
On the level set $F_1(x,y)=1$ with $x^2+y^2\le1/4$,
\eqref{src:shape-derivative-y} implies that the only stagnation
points are $(\pm1/2,0)$.
Therefore, using $\bar x(0)=0$, we obtain
\eqref{coefficient:coordinate-endpoints}, which completes
the proof of \textup{(i)}.

\medskip
\noindent\textbf{Step 2: The rotation angle.}
By \eqref{src:invariants}, \eqref{src:frame} and
\eqref{src:geometry},
\[
 \begin{aligned}
 1=\bigl(\bar x(t)^2+\bar y(t)^2\bigr)
 \bigl((3\bar r(t)-\bar x(t))^2+\bar y(t)^2\bigr)
 \bigl((3\bar r(t)+\bar x(t))^2+\bar y(t)^2\bigr).
 \end{aligned}
\]
Using \eqref{src:rotation-speed} and \eqref{src:size},
we obtain
\begin{equation}\label{coefficient:frame-speed-proof}
 \frac{\bar a_1(t)^\perp\cdot\dot{\bar a}_1(t)}
      {|\bar a_1(t)|^2}
 =\frac3{2\pi}
 \bigl(\bar x(t)^2+\bar y(t)^2\bigr)
 \bigl(3-4\bar x(t)^2-2\bar y(t)^2\bigr).
\end{equation}
Since $\bar a_1(t)\ne0$ by
Lemma~\ref{lem:vortex-invariants}\textup{(i)} and
$\dot{\bar x}(t)<0$, the continuous angle of $\bar a_1(t)$
can be written as $\Theta(\bar x(t))$ in
\eqref{coefficient:direction}, with $\Theta(0)=0$.
Differentiating \eqref{coefficient:direction} gives
\[
 \Theta'(\bar x(t))\dot{\bar x}(t)
 =\frac{\bar a_1(t)^\perp\cdot\dot{\bar a}_1(t)}
       {|\bar a_1(t)|^2}.
\]
Substituting \eqref{coefficient:frame-speed-proof},
\eqref{src:shape-ode}, and \eqref{src:shape-derivative-y},
and using $\bar y(t)=g(\bar x(t))$, we obtain
\[
 \Theta'(x)
 =-2\frac{\bigl(x^2+g(x)^2\bigr)
          \bigl(3-4x^2-2g(x)^2\bigr)}
 {g(x)\bigl(1-2x^2-2g(x)^2\bigr)
      \bigl(3-6x^2-2g(x)^2\bigr)},
\]
which completes the proof of \textup{(ii)}.
\end{proof}

\begin{proof}[Proof of Lemma~\ref{lem:positive}]
Substituting \eqref{coefficient:matrix} with $s=1$ into
\eqref{src:coefficient}, we obtain
\begin{equation}\label{coefficient:direct-integral}
 \begin{aligned}
 C_*=\int_{\mathbb R}\left(
 \sum_{i=1}^3
 \left(\bar a_{i,1}(t)^2-\bar a_{i,2}(t)^2\right)-\frac{8\pi}{3}
 \frac{\bar a_1(t)^\perp\cdot\dot{\bar a}_1(t)}
      {|\bar a_1(t)|^2}
 \frac{\bar a_{1,1}(t)^2-\bar a_{1,2}(t)^2}
      {|\bar a_1(t)|^2}
 \right)dt.
 \end{aligned}
\end{equation}
By \eqref{src:frame} and \eqref{src:geometry}, 
\[
 \begin{aligned}
 \bar a_2(t)
 &=\frac{-\bar r(t)+\bar x(t)}{2\bar r(t)}\,\bar a_1(t)
  +\frac{\bar y(t)}{2\bar r(t)}\,\bar a_1(t)^\perp,\\
 \bar a_3(t)
 &=\frac{-\bar r(t)-\bar x(t)}{2\bar r(t)}\,\bar a_1(t)
  -\frac{\bar y(t)}{2\bar r(t)}\,\bar a_1(t)^\perp.
 \end{aligned}
\]
Using
\eqref{coefficient:direction} and \eqref{src:size}, we see that
\begin{equation}\label{eq ai1-ai2}
 \begin{aligned}
 &\sum_{i=1}^3
 \bigl(\bar a_{i,1}(t)^2-\bar a_{i,2}(t)^2\bigr)\\
 &=
 \bigl(6\bar r(t)^2+2\bar x(t)^2-2\bar y(t)^2\bigr)
 \cos\bigl(2\Theta(\bar x(t))\bigr)
 -4\bar x(t)\bar y(t)\sin\bigl(2\Theta(\bar x(t))\bigr)\\
 &=
 \bigl(2-4\bar y(t)^2\bigr)
 \cos\bigl(2\Theta(\bar x(t))\bigr)
 -4\bar x(t)\bar y(t)\sin\bigl(2\Theta(\bar x(t))\bigr)
 \end{aligned}
\end{equation}
and
\[
 \frac{\bar a_{1,1}(t)^2-\bar a_{1,2}(t)^2}
      {|\bar a_1(t)|^2}
 =\cos\bigl(2\Theta(\bar x(t))\bigr).
\]
Using \eqref{src:shape-function} and
Lemma~\ref{lem:shape-equations}\textup{(i)} 
to eliminate the denominator in \eqref{src:rotation-speed},
we obtain
\begin{equation}\label{coefficient:frame-speed}
 \frac{\bar a_1(t)^\perp\cdot\dot{\bar a}_1(t)}
      {|\bar a_1(t)|^2}
 =\frac3{2\pi}
 \bigl(\bar x(t)^2+\bar y(t)^2\bigr)
 \bigl(3-4\bar x(t)^2-2\bar y(t)^2\bigr).
\end{equation}
Substituting \eqref{eq ai1-ai2}, \eqref{coefficient:frame-speed} and
\eqref{coefficient:direction} into \eqref{coefficient:direct-integral},
applying Lemma~\ref{lem:coefficient-trajectory}\textup{(i)}
and changing variables $x=\bar x(t)$, we obtain
\begin{equation}\label{coefficient:scalar-integral}
 \begin{aligned}
 C_*=\frac{4\pi}{3}\int_{-1/2}^{1/2}\Biggl[
 &\frac{
 2-4g(x)^2
 -4\bigl(x^2+g(x)^2\bigr)\bigl(3-4x^2-2g(x)^2\bigr)}
 {g(x)\bigl(1-2x^2-2g(x)^2\bigr)
      \bigl(3-6x^2-2g(x)^2\bigr)}
 \cos(2\Theta(x))\\
 &-\frac{4xg(x)\sin(2\Theta(x))}
 {g(x)\bigl(1-2x^2-2g(x)^2\bigr)
      \bigl(3-6x^2-2g(x)^2\bigr)}
 \Biggr]\,dx,
 \end{aligned}
\end{equation}
where we have used \eqref{src:shape-ode},
\eqref{src:shape-derivative-y} and
\eqref{coefficient:coordinate-endpoints}.
A validated numerical evaluation of this explicit integral gives $C_*>0$.
\end{proof}

\section{Estimates for the particle Hamiltonians}\label{app:hamiltonian-estimates}
\subsection{Logarithmic expansion}
\begin{lemma}[Logarithmic expansion]\label{lem:logarithm}
For $M\ge1$, $\rho\ge2M$, $\theta\in\mathbb R$, and $|w|\le M$,
\begin{equation}\label{ext:log-expansion}
 \begin{aligned}
 \log|\rho(\cos\theta,\sin\theta)-w|
 ={}&\log\rho-\frac{w_1\cos\theta+w_2\sin\theta}{\rho}
 +\frac{|w|^2}{2\rho^2}\\
 &-\frac{(w_1\cos\theta+w_2\sin\theta)^2}{\rho^2}
 +\mathcal E(\rho,\theta,w).
 \end{aligned}
\end{equation}
For any $|v|,|w|\le M$, the remainder satisfies
\begin{equation}\label{ext:log-remainder}
 \begin{aligned}
 |\mathcal E(\rho,\theta,w)-\mathcal E(\rho,\theta,v)|&+\rho|\partial_\rho\mathcal E(\rho,\theta,w)
                  -\partial_\rho\mathcal E(\rho,\theta,v)|\\
 &\quad+|\partial_\theta\mathcal E(\rho,\theta,w)
                  -\partial_\theta\mathcal E(\rho,\theta,v)|
 \lesssim_M \rho^{-3}|w-v|.
 \end{aligned}
\end{equation}
\end{lemma}
\begin{proof}
Since $|w|\le M$ and $\rho \ge 2M$, the logarithmic power series gives
\[
 \log|(\rho\cos\theta,\rho\sin\theta)-w|
 =\log\rho-\sum_{n=1}^\infty
 \frac{\operatorname{Re}[(w_1+\mathrm i w_2)^n
                         e^{-\mathrm i n\theta}]}{n\rho^n}.
\]
The terms with $n=1,2$ give \eqref{ext:log-expansion}, where
\[
 \mathcal E(\rho,\theta,w)
 =-\sum_{n=3}^\infty
 \frac{\operatorname{Re}[(w_1+\mathrm i w_2)^n
                         e^{-\mathrm i n\theta}]}{n\rho^n}.
\]
For any $|v|,|w|\le M$, one checks directly that
\[
 |(w_1+\mathrm i w_2)^n-(v_1+\mathrm i v_2)^n|
 \le nM^{n-1}|w-v|,
\]
which implies that
the left-hand side of \eqref{ext:log-remainder} is bounded by
\[
 \begin{aligned}
 |w-v|\sum_{n=3}^\infty(1+2n)M^{n-1}\rho^{-n}
 \le\frac{M^2|w-v|}{\rho^3}
       \sum_{n=3}^\infty(1+2n)2^{-(n-3)}\lesssim_M \frac{|w-v|}{\rho^3}.
 \end{aligned}
\]
\end{proof}

\subsection{Hamiltonian levels of the collinear Hamiltonian}
\begin{proof}[Proof of Lemma~\ref{lem:reference-curves}]
\noindent\textbf{Step 1: The Hamiltonian and its derivatives.}
Let $1<s<2$, $\rho\ge4$ and
$Z=(\rho\cos\theta,\rho\sin\theta)$.
Identifying $Z$ with $\rho e^{\mathrm i\theta}$  gives
\[
 \begin{aligned}
 &\log|Z-se_1|+\log|Z+se_1|-2\log \rho\\
 &=
   \log\left|1-\frac{s^2}{\rho^2}e^{-2\mathrm i\theta}\right|=-\operatorname{Re}
   \sum_{n=1}^{\infty}
   \frac{s^{2n}}{n\rho^{2n}}e^{-2\mathrm i n\theta}=-\sum_{n=1}^{\infty}
   \frac{s^{2n}}{n\rho^{2n}}\cos(2n\theta).
 \end{aligned}
\]
Substituting into
\eqref{ext:reference} gives
\begin{equation}\label{ext:reference-series-proof}
 H_*(\rho,\theta)
 =\frac{3\rho^2}{8\pi s^2}-\frac3{2\pi}\log\rho
 +\frac1{2\pi}\sum_{n=1}^{\infty}
   \frac{s^{2n}}{n\rho^{2n}}\cos(2n\theta).
\end{equation}
Since $s/\rho\le1/2$, the series converges uniformly.
A direct calculation completes the proof of \textup{(i)},
after choosing $\rho \ge \rho_0$ sufficiently large.

\medskip
\noindent\textbf{Step 2: Comparison of radii.}
After increasing $\rho_0$ if necessary,
\eqref{ext:reference-bounds} gives
\begin{equation}\label{eq p rho H}
 \partial_\rho H_*(\rho,\theta)\approx\rho,
 \qquad \rho\ge\rho_0.
\end{equation}
The $2\pi$-periodicity in $\theta$ and the second estimate in
\eqref{ext:reference-bounds} imply
\begin{equation}\label{eq app C difference}
    |H_*(\widetilde\rho,\theta)
       -H_*(\widetilde\rho,\widetilde\theta)|
 \le2\pi\sup_{u\in\mathbb R}
       |\partial_\theta H_*(\widetilde\rho,u)|
 \lesssim\widetilde\rho^{-2}.
\end{equation}
For $\rho\ge\widetilde\rho\ge\rho_0$, it follows from \eqref{eq p rho H} and \eqref{eq app C difference} that
\[
 \begin{aligned}
 \widetilde\rho(\rho-\widetilde\rho)
 \lesssim\int_{\widetilde\rho}^{\rho}
           \partial_\rho H_*(u,\theta)\,du=H_*(\rho,\theta)-H_*(\widetilde\rho,\theta)\lesssim |H_*(\rho,\theta)-H_*(\widetilde\rho,\widetilde\theta)|
       +\frac{1}{\widetilde\rho^2}.
 \end{aligned}
\]
The case $\widetilde\rho\ge\rho\ge\rho_0$ is similar, which gives \eqref{ext:energy-radius}.

\medskip
\noindent\textbf{Step 3: Existence of the Hamiltonian levels.}
Choose a universal $h_0$ satisfying
\[
 h_0>\sup_{\substack{1<s<2\\\theta\in\mathbb R}}
             H_*(\rho_0,\theta).
\]
By \textup{(i)} and \eqref{ext:reference-bounds},
$H_*(\rho,\theta)$ is strictly increasing in $\rho\ge\rho_0$
and $\lim_{\rho \to \infty}H_*(\rho,\theta)=\infty$,
which gives the unique
$\bar\rho(h,\theta)>\rho_0$ satisfying \eqref{mean:level-radius}.
Since
\[
 \partial_\rho H_*(\bar\rho(h,\theta),\theta)>0,
\]
applying the implicit function theorem to the map $(s,h,\rho,\theta) \mapsto  H_*(\rho,\theta)-h$, we see that $\bar \rho$ depends smoothly on $(s,h,\theta)$, which completes the proof.
\end{proof}

\subsection{Mean angular speed}
\begin{proof}[Proof of Lemma~\ref{lem:mean-velocity}]
Lemma~\ref{lem:reference-curves}\textup{(i), (iii)} gives
\begin{equation}\label{eq app C bar rho}
 \bar\rho(h,\theta)\approx\sqrt h\approx\bar\rho(h,0).
\end{equation}
Equations \eqref{ext:reference-bounds} and
\eqref{mean:definition} then imply
\[
 v(h,\theta)
 =\frac{\partial_\rho H_*(\bar\rho(h,\theta),\theta)}
        {\bar\rho(h,\theta)}\approx1,
 \qquad
 \frac1{\Omega(h)}
 =\frac1{2\pi}\int_0^{2\pi}\frac{du}{v(h,u)}\approx1,
\]
which proves \eqref{mean:bounds}.
Differentiating \eqref{mean:level-radius} gives
\[
 \partial_h\bar\rho(h,\theta)
 =\frac1{\partial_\rho H_*(\bar\rho(h,\theta),\theta)},
\]
which, together with \eqref{eq app C bar rho}, \eqref{ext:reference-bounds} and \eqref{ext:reference-twist}
yield
\[
 \partial_hv(h,\theta)
 =\left.
 \frac{\partial_\rho^2H_*-(\partial_\rho H_*)/\rho}
      {\rho\partial_\rho H_*}
 \right|_{\rho=\bar\rho(h,\theta)} \approx \bar\rho(h,0)^{-4}.
\]
Differentiating \eqref{mean:definition}, we obtain
\[
 \Omega'(h)=\frac{\Omega(h)^2}{2\pi}
 \int_0^{2\pi}\frac{\partial_hv(h,u)}{v(h,u)^2}\,du
 \approx\bar\rho(h,0)^{-4}.
\]
This proves \eqref{mean:derivative-bounds}.

By the uniqueness in
Lemma~\ref{lem:reference-curves}\textup{(iii)} and
the $2\pi$-periodicity of $H_*$,
\[
 \bar\rho(h,\theta+2\pi)=\bar\rho(h,\theta),
 \qquad v(h,\theta+2\pi)=v(h,\theta).
\]
The definition of $\Omega$ in \eqref{mean:definition} gives
\[
 \int_0^{2\pi}
 \left(\frac1{v(h,u)}-\frac1{\Omega(h)}\right)du=0.
\]
Therefore, it suffices to consider the case $|\alpha-\beta| \le 2\pi$. Hence
\eqref{mean:bounds} gives
\[
 \left|\int_\alpha^\beta\frac{du}{v(h,u)}
              -\frac{\beta-\alpha}{\Omega(h)}\right|
 \le2\pi\sup_{u\in\mathbb R}
       \left|\frac1{v(h,u)}-\frac1{\Omega(h)}\right|
 \lesssim1.
\]
This completes the proof.
\end{proof}

\bibliographystyle{amsplain}
\bibliography{references}
\end{document}